\documentclass[12pt]{article}
\usepackage{mathrsfs}
\usepackage{amsmath}
\usepackage{amsfonts}
\usepackage{amsthm}
\usepackage{graphicx}
\usepackage{float}
\usepackage{tikz}
\usepackage{pgf}
\usepackage{color}
\usepackage{geometry}                
\usepackage{graphicx}
\usepackage{amssymb}
\newcommand{\R}{\mathbb{R}}
\newcommand{\Z}{\mathbb{Z}}
\newcommand{\N}{\mathbb{N}}
\newcommand{\eps}{\varepsilon}

\newcounter{number}
\newtheorem*{thmA}{Theorem A}
\newtheorem*{thmB}{Theorem B}
\newtheorem*{thmB'}{Theorem B$'$}
\newtheorem*{thmC}{Theorem C}

\newtheorem{theorem}{Theorem}
\newtheorem{lemma}{Lemma}
\newtheorem{proposition}{Proposition}
\newtheorem{corollary}{Corollary}
\newtheorem{remark}{Remark}
\newtheorem{definition}{Definition}
\newcommand{\begle}{\begin{lemma}\hspace{1em}}
\newcommand{\ele}{\end{lemma}
                       \addtocounter{number}{1}}
\newcommand{\begde}{\begin{definition}\hspace{1em}}
\newcommand{\ede}{\end{definition}
                           \addtocounter{number}{1}}
\newcommand{\begth}{\begin{theorem}\hspace{1em}}
\newcommand{\eeth}{\end{theorem}
   \addtocounter{number}{1}}

\newcommand{\begpro}{\begin{proposition}\hspace{1em}}
\newcommand{\epro}{\end{proposition}
   \addtocounter{number}{1}}
\newcommand{\begrem}{\begin{remark}\hspace{1em}}
\newcommand{\erem}{\end{remark}
   \addtocounter{number}{1}}

\newcommand{\begco}{\begin{corollary}\hspace{1em}}
\newcommand{\eco}{\end{corollary}
   \addtocounter{number}{1}}
\newcommand{\begeq}{\begin{equation}}
\newcommand{\eeq}{\end{equation}}
\renewcommand{\thelemma}{\thesection.\thenumber}

\newcommand{\newsection}[1]{    \section{#1}
                               \setcounter{equation}{0}
                                \setcounter{number}{1} }
\newcommand{\emark}{\hfill\hspace{1em}\framebox[0.5em]{\rule{0cm}{0em}}
\vspace{2mm}}

\renewcommand{\theequation}{\arabic{section}.\arabic{equation}}

\newcommand{\refer}[1]{~(\ref{#1})}
\newcommand{\wbf}{{\mathbf{w}}}

\newcommand{\ubf}{{\mathbf{u}}}
\newcommand{\zbf}{{\mathbf{z}}}
\newcommand{\xbf}{{\mathbf{x}}}
\newcommand{\ybf}{{\mathbf{y}}}

\newcommand{\one}{{\mathbf{1}}}
\newcommand{\zero}{{\mathbf{0}}}

\newcommand{\mcalX}{\mathcal{X}}
\newcommand{\nubar}{{\bar{\nu}}}

\newcommand{\pf}{{\bf Proof:}\hspace{1em}}

\usepackage{soul}

\title{Depinning of discommensurations for tilted Frenkel-Kontorova chains}
\author{
{\small $\mbox{C. Baesens}^1$ \quad $\mbox{R. S. MacKay}^1$\quad $\mbox{W.-X. Qin}^2$
}
 \\
 {\small 1.  Mathematics Institute,  University of Warwick, Coventry CV4 7AL, UK\ \ \  }\\
{\small 2. Department of Mathematics, Soochow University, Suzhou, 215006,
 China}
}
\date{}                                         

\begin{document}
\maketitle
\begin{abstract}
For an untilted Frenkel-Kontorova chain and any rational $p/q$, Aubry and Mather proved there are minimising equilibrium states that are left- and right-asymptotic to neighbouring pairs of spatially periodic minimisers of type $(p,q)$.  They are known as {\em discommensurations} (or kinks or fronts), {\em advancing }if the right-asymptotic equilibrium is to the right of the left-asymptotic one, {\em retreating} otherwise.  Following work of Middleton, Floria \& Mazo and Baesens \& MacKay, there is a threshold tilt $F_d(p/q)\ge 0$ up to which there continue to be periodic equilibria of type $(p,q)$ and above which there 
is a globally attracting periodically sliding solution in the space of sequences of type $(p,q)$.

In this paper, we prove that there are values $F_d(p/q\pm)$ of tilt with $0\le F_d(p/q\pm) \le F_d(p/q)$, generically positive and less than $F_d(p/q)$, up to which there continue to be equilibrium advancing or retreating discommensurations, respectively,
and such that for $F_d(p/q\pm) < F < F_d(p/q)$ there
are periodically sliding discommensurations, apart perhaps from exceptional cases with both a degenerate type $(p,q)$ equilibrium and a degenerate advancing equilibrium discommensuration.  We give examples, however, to show that equilibrium and periodically sliding discommensurations may co-exist, both above and below $F_d(p/q\pm)$,
so the case of discommensurations is not as clean as that of periodic configurations.

On the way, we prove that $F_d(\omega) \to F_d(p/q\pm)$ as $\omega \searrow p/q$ or $\nearrow p/q$ respectively.  Finally, we prove that $F_d(p/q\pm)=0$ is equivalent to the existence of a rotational invariant circle consisting of periodic orbits of type $(p,q)$ and right-going (respectively left-going) separatrices, for the corresponding twist map on the cylinder.
\end{abstract}

{\noindent \bf Keywords:}\hspace{0.5em}  Discommensuration,  Depinning Transition,  Aubry-Mather Theory, Frenkel-Kontorova Chain, Invariant Circle.


\newsection{Introduction}\label{sec:int}

Frenkel-Kontorova (FK) chains are models for one-dimensional structures in solid-state physics, like charge density waves~\cite{Gr,Th}.  They are defined by a nearest-neighbour energy function $h: \R\times\R \to \R$ with periodicity property
\begin{equation}
h(x+1,x'+1) = h(x,x'),
\label{eq:per}
\end{equation}
which is $C^2$ and has negative mixed second derivative $h_{12}(x,x') \le -c < 0$ (weakened conditions can be considered \cite{Ban}).
A standard example is
\begin{equation}
h(x,x') = \frac12 (x'-x)^2 + \frac{k}{4\pi^2} \cos{2\pi x}.
\label{eq:stdFK}
\end{equation}
The formal energy of an infinite chain with state $\xbf = (x_n)_{n \in \Z} \in \R^\Z$ is $$W(\xbf) = \sum_{n\in \Z} h(x_n, x_{n+1}),$$
though in general the sum does not converge.
Equilibria of an FK chain are the states $\xbf=(x_n)$ for which $\frac{\partial W}{\partial x_n} = 0$ for all $n \in \Z$, i.e.~$$h_2(x_{n-1},x_n)+h_1(x_n,x_{n+1})=0,$$
where subscript $i$ on $h$ denotes the partial derivative with respect to the $i^{th}$ argument.
The gradient vector field of an FK chain is
$$\dot{x}_n = -\frac{\partial W}{\partial x_n} = -h_2(x_{n-1},x_n)-h_1(x_n,x_{n+1}).$$

Tilted FK chains are defined by a nearest-neighbour energy function $h_F$ in $C^2$ with negative mixed second derivative, but the periodicity property is replaced by $$h_F(x+1,x'+1)-h_F(x,x')=-F,$$ where $F$ is called the {\em tilt}.
Defining $h(x,x')=h_F(x,x')+Fx$, we see that $h$ satisfies the periodicity property\refer{eq:per}.
So we will consider families of the form
\begin{equation}
h_F(x,x')=h(x,x')-Fx,
\label{eq:hF}
\end{equation}
with $h$ satisfying the periodicity condition\refer{eq:per}.  Without loss of generality we take $F\ge 0$ (else reverse the signs of $x,x'$).
Thus the gradient vector field of a tilted FK chain is
\begin{equation}
\dot{x}_n = -h_2 (x_{n-1},x_n)-h_1(x_n,x_{n+1})+F.
\label{tilt}
\end{equation}
It induces a flow $\phi_F$ on suitable spaces of sequences (see Sec.~\ref{sec:gradflow}).
Equilibria $\xbf=(x_n)$ satisfy
\begin{equation}\label{=F}
h_2 (x_{n-1},x_n)+h_1(x_n,x_{n+1})=F,
\end{equation}
corresponding to orbits of an  area-preserving monotone twist map $\Phi_F$: $(x_n,p_n) \mapsto (x_{n+1},p_{n+1})$ on the cylinder, determined implicitly by
\begin{equation}\label{eq:apmap}
\left\{
\begin{array}{lll}
p_n &=& -h_1(x_n,x_{n+1})+F, \\
p_{n+1}&=& h_2(x_n,x_{n+1}).
\end{array}
\right.
\end{equation}

Endow the space $\R^\Z$ of sequences with the partial order $\xbf \le \ybf$ if $x_n \le y_n$ for all $n \in \Z$.
We write $\xbf < \ybf$ if $\xbf \le \ybf$ and $\xbf \ne \ybf$ (i.e.~there is at least one $n$ for which $x_n < y_n$).
We write $\xbf \ll \ybf$ if $x_n < y_n$ for all $n \in \Z$, and say $\xbf$ is {\em strictly} less than $\ybf$.

The translations $T_{qp}$ in $\R^{\Z}$ are defined for $p,q\in\Z$ by
$$(T_{qp}\xbf)_n=x_{n-q}+p .$$
A configuration $\xbf=(x_n)$ is said to be {\em Birkhoff} (or {\em weakly ordered}) if the set of translates $\{T_{qp}\xbf\,|\,p, q \in \Z\}$ is totally ordered with respect to the partial order defined above on $\R^\Z$.  In particular this implies that each Birkhoff configuration $\xbf=(x_n)$ has a {\em mean spacing} $\omega = \lim_{n\to\pm \infty} (x_n-x_0)/n \in \R$ and that $x_n-n\omega$ is bounded~\cite{BM}.

For $p\in \Z, q \in \N$, we say a configuration $\xbf$ is of {\em type $(p,q)$} if $T_{qp}\xbf=\xbf$. So it is periodic modulo translation by $\Z$.  Not every Birkhoff configuration of mean spacing $p/q$ is periodic, however.  There can be ``discommensurations'' of mean spacing $p/q$ with $T_{qp}\xbf<\xbf$ or $T_{qp}\xbf>\xbf$, called advancing and retreating respectively (but the names refer to spatial dependence, not time-dependence).  Discommensurations are left- and right-asymptotic to type $(p,q)$ configurations.

For each $\omega \in \R$ irrational, or $\omega=p/q$ rational (which we take in lowest terms and with $q>0$), there is a {\em depinning force} $F_d(\omega) \ge 0$  for the tilted FK chain\refer{tilt} such that for $0\le F \le F_d(\omega)$ there is a Birkhoff equilibrium of mean spacing $\omega$ (of type $(p,q)$ if $\omega$ is rational $p/q$), and for $F>F_d(\omega)$ there are no Birkhoff equilibria with mean spacing $\omega$.  Moreover, for $F>F_d(\omega)$ there is a time-periodic solution $\xbf(t)=(x_n(t))_{n\in\Z}$ of mean spacing $\omega$ (type $(p,q)$ if $\omega$ is rational $p/q$) of the gradient flow in the sense that $x_n(t+T) = x_n(t)+1$ for some $T=T(\omega)>0$;
 $v(\omega)=1/T$ is called the {\em average velocity}
and one can write
$$x_n(t) = X(n\omega+vt+\alpha)$$
for a ``dynamical hull function'' $X$ which is increasing, differentiable and diagonally periodic, $X(\alpha+1)=X(\alpha)+1$, $\forall\,\alpha\in\R$.
There are no such solutions for $F\in [0,F_d(\omega)]$.
These results were formulated in \cite{FM} (following \cite{Mi} for the case of finite chains).  They were proved in \cite{BM} for the rational case.  The claimed proof in \cite{BM} for the irrational case contains an error, pointed out in \cite{BM2}, but the result was proved in \cite{Q} for the special example\refer{eq:stdFK}, and in \cite{QW1} for general $h$ with bounded second derivative, and their proof can be extended to the case without bounded second derivative by the remark at the beginning of our Section~2.1.

In the present paper, we study what happens to the discommensurations whose existence was proved by Aubry and Mather for $F=0$, when a tilt is added.
We generalise the definition of discommensuration to Birkhoff configurations left- and right-asymptotic to different Birkhoff periodic states of the same type $(p,q)$ (i.e.~we do not require minimisers).
Numerics have been performed on the question by several groups, e.g.~\cite{FM}, and some conjectures were given in \cite{BM}.
Some results have been proved in \cite{QW1,QW2},
and we adapt their approach.

We always assume in this paper that the energy function $h: \R\times\R \to \R$ with periodicity condition\refer{eq:per} is $C^2$ and has negative mixed second derivative $h_{12}(x,x') \le -c < 0$.

Our main conclusions are the following Theorems A, B and C.

\begin{thmA}\label{thmA}
The limits $ \lim_{\omega \searrow p/q} F_d(\omega) =F_d(p/q+)$ and $\lim_{\omega \nearrow p/q}F_d(\omega)=F_d(p/q-)$ exist. Moreover, $0\le F_d(p/q\pm) \le F_d(p/q)$, and the inequalities are generically strict.
\end{thmA}

Note that the continuity properties of $F_d$ can be summarised in terms of Mather's concept of the space $\mathcal{R}$ of rotation symbols~\cite{Mat3}.  The rotation symbols are obtained from the reals by blowing up each rational $p/q$ into three points $p/q-, p/q, p/q+$, with order $\omega_1 < p/q- < p/q < p/q+ < \omega_2$ for all reals $\omega_1 < p/q < \omega_2$.  The topology on $\mathcal{R}$ is that generated by the open intervals in this order.  The rationals become isolated points in the gaps of a $\Z$-cover of a Cantor set.  Mather proved that his $\Delta W$ function is continuous from $\mathcal{R}$ to $\R$~\cite{Mat2}.  Similarly, we prove that $F_d: \mathcal{R} \to \R$ is continuous.

\begco\label{co:conti}
The depinning force $F_d:{\cal R}\to\R$ is continuous.
\eco



Assume that $\xbf^-<\xbf<\xbf^+$ are Birkhoff equilibria of\refer{tilt} and both $\xbf^-$ and $\xbf^+$ are $(p,q)$-periodic. If
$|x_n-x^-_n| \to 0$ as $n \to -\infty$ and $|x_n-x^+_n| \to 0$ as $n \to +\infty$, then $\xbf$ is called an {\em equilibrium advancing discommensuration} of type $p/q+$; similarly if $|x_n-x^+_n| \to 0$ as $n \to -\infty$ and $|x_n-x^-_n| \to 0$ as $n \to +\infty$, it is called an {\em equilibrium retreating discommensuration} of type $p/q-$.
They are heteroclinic orbits of the associated map $\Phi_F$, between the $(p,q)$-periodic orbits.  Discommensurations can alternatively be described as defects, the advancing case corresponding to a vacancy and the retreating case to an interstitial.

The definitions can be extended to time-dependent solutions $\xbf(t)$ of\refer{tilt} as follows.  Given Birkhoff equilibria $\xbf^- < \xbf^+$ of type $(p,q)$, $\xbf(t)$ is an {\em advancing discommensuration} if $\xbf^- < \xbf(t) < \xbf^+$ for all $t\in \R$ and $|x_n(t)-x_n^\pm|\to 0$ as $n\to \pm \infty$ respectively.  It is a {\em retreating discommensuration} if $\xbf^- < \xbf(t) < \xbf^+$ for all $t\in \R$ and $|x_n(t)-x_n^\pm|\to 0$ as $n\to \mp \infty$ respectively.  It is {\em periodically sliding} if there exists $T>0$ such that $\xbf(t+T) = T_{-q,-p}\xbf(t)$ in the advancing case, $T_{qp}\xbf(t)$ in the retreating case.

The significance of the quantities $F_d(p/q\pm)$ obtained in Theorem A is that 
for $F\le F_d(p/q\pm)$ there are equilibrium advancing or retreating discommensurations, respectively,  and  for 
$F_d(p/q\pm)<F<F_d(p/q))$ there are periodically sliding discommensurations,
apart perhaps from an exceptional case, as we now state.

\begin{thmB}\label{thmB}
 For $0< F \le F_d(p/q+)$, there exist for\refer{tilt} equilibrium advancing discommensurations of type $p/q+$. 

  For $F_d(p/q+) < F < F_d(p/q)$, there are periodically sliding advancing discommensurations or there is both a degenerate equilibrium of type $(p,q)$ and a degenerate advancing equilibrium discommensuration of type $p/q+$.

 The analogous results hold for $p/q-$.
\end{thmB}

Note that if in particular $p=0$ and $q=1$, then a periodically sliding advancing discommensuration 
is in the usual sense a travelling front solution and an equilibrium advancing discommensuration is a stationary front. 

The condition that there be both a degenerate equilibrium of type $(p,q)$ and a degenerate advancing equilibrium discommensuration is codimension-2,  
 thus can be considered unlikely in a family parametrised by one parameter, such as the tilt $F$ here.
Perhaps one could strengthen Theorem B further by reducing this to an even higher codimension situation. 

It is important to note that
Theorem B does not say that there are no equilibrium advancing discommensurations for $F_d(p/q+)<F<F_d(p/q)$, nor that there are no periodically sliding ones for $0<F<F_d(p/q+)$.  Indeed in Appendix B we make {examples} with coexistence of equilibrium and periodically sliding ones of rotation number $0/1+$, {for $F>F_d(0/1+)$ and for $F<F_d(0/1+)$.}

\vskip 1ex
{\noindent\bf Remark:} We compare our approach and results with those represented by \cite{CMPS, MP, BCC}.  They study the gradient dynamics of chains with a bistable potential at each site.  When the height difference between the wells is small enough (or equivalently, the coupling along the chain is weak enough) they obtain stationary fronts; when it is large enough they obtain travelling fronts.  These are analogous to our discommensurations and the transition at $F_d(0/1+)$.
In contrast to us, they require non-degeneracy of the minima of the potential and fairly heavy analysis.  We use their results, however, in Appendix B.

\vspace{2mm}

For irrational $\omega$ or rational $\omega=p/q$, it is shown~\cite{QW1} that $F_d(\omega)=0$ if and only if there exists a rotational invariant circle with rotation number $\omega$ (consisting of $(p,q)$-periodic orbits in the rational case) for the corresponding monotone twist map $\Phi_0$ ($F=0$ in\refer{eq:apmap}). Similarly, we can also give via $F_d(p/q\pm)$ a criterion for the existence of invariant circles consisting of $(p,q)$-periodic orbits and separatrices.

\begin{thmC}\label{thmC}
The twist map\refer{eq:apmap} with $F=0$ has a rotational invariant circle on the cylinder consisting of $(p,q)$-periodic orbits, or of $(p,q)$-periodic orbits together with orbits corresponding to equilibrium advancing (retreating) discommensurations, if and only if $F_d(p/q+)=0$ ($F_d(p/q-)=0$).
\end{thmC}

Note that if the periodic orbits are hyperbolic (as is generic) the orbit $(x_n,p_n)_{n\in\Z}$ of $\Phi_0$ corresponding to the equilibrium advancing discommensuration $\xbf=(x_n)$ lies on the unstable manifold of the periodic orbit corresponding to $\xbf^-$ and stable manifold of the  orbit corresponding to $\xbf^+$, where $\xbf^-$ and $\xbf^+$ are respectively the left- and right-asymptotic equilibria of  $\xbf$. Therefore, in this case the invariant circle consists of $(p,q)$-periodic orbits and separatrices.


Theorem C does not totally resolve the question of existence of rational rotational invariant circles for $\Phi_0$, because there could be rotational invariant circles formed of branches of invariant manifold of periodic orbits in opposite directions (see  Appendix D).

The paper \cite{QW1} constructed (whatever the tilt $F$) a totally ordered closed set of states homeomorphic to $\R$, invariant under both the gradient flow and the translations.  It is called an {\em invariant ordered circle} (IOC), because it is totally ordered, its quotient by $T_{01}$ is a circle, and it is invariant under all the translations and the gradient flow.  The method is a generalisation from the untilted case of Hall's ridge curves, which gave an alternative construction of the results of Aubry and Mather, were developed further under the name of ghost circles by \cite{Go} and used by \cite{LC}.  The key result of \cite{QW1} is that in the space of sequences of type $(p,q)$ there is an IOC, using Schauder's fixed point theorem for the gradient flow on the subset of totally ordered curves.  By taking a limit point as $p/q \to \omega$ irrational, an IOC is obtained for mean spacing $\omega$.

In \cite{QW2}, the result was extended under the assumption of continuity of $F_d$ at $p/q$ to construct equilibria of types $p/q\pm$ for $F<F_d(p/q)$, meaning states with $T_{qp}\xbf\le \xbf$, respectively $\ge \xbf$, and of mean spacing $p/q$.
This assumption is very restrictive, however, as Theorem A shows that generically $F_d$ is not continuous at rationals. 

Here we prove our main results by constructing IOCs of types $p/q\pm$ (to be defined in Section~\ref{sec:p/q+}) for $0\leq F< F_d(p/q)$.

These results partially answer question (iii) of \cite{BM} (and Q.4 of \cite{B}) and replace the paper of Baesens \& Floria promised there.

Our paper starts by recalling some results on the gradient flow and depinning force, the classic Aubry-Mather theory for untilted FK chains (Section~\ref{sec:pre}), and the definition and theory of IOCs, together with the construction of IOCs of type $p/q\pm$ (Section~\ref{sec:p/q+}).  Then we present the proofs of our main conclusions (Section~\ref{sec:pf}). 
 We discuss the significance of the results in Section~\ref{sec:dis}, and end with several appendices. 



\newsection{Preliminaries}\label{sec:pre}
\subsection{Gradient flow and depinning force} \label{sec:gradflow}
In this paper we always assume that $h$ is $C^2$ and has negative mixed second derivative $h_{12}(x,x') \le -c < 0$.
We do not require $h$ to have bounded second derivative, but given any integers $M<N$ one can modify $h(x,x')$ for $x'-x<M$ and for $x'-x>N$ to make the second derivative bounded, while keeping $h\in C^2$ and $h_{12}\le -c$.  An explicit choice is given in  Appendix E.  Our interest will be in sequences $\xbf$ satisfying $M\le x_{n+1}-x_n\le N$, so for them the modification is zero.

System\refer{tilt} defines a flow on various spaces.  One set of examples is the spaces $\mathcal{X}_\omega$, $\omega \in \R$, of sequences $\xbf=(x_n)$ with $x_n - n\omega$ bounded and $\ell_\infty$ (supremum) metric.  Then the gradient vector field is $C^1$ on $\mathcal{X}_\omega$ and induces a $C^1$ flow on $\mathcal{X}_\omega$.  We will need, however, to compare the flow for different $\omega$ in product topology.  For this purpose we follow \cite{Go2,Q,QW1} by using the space $\mathcal{X}$ of sequences $\xbf=(x_n)\in\R^{\Z}$ with norm
$$
\|\xbf\|=\sum_{n\in \Z} 2^{-|n|} |x_n| < \infty,
$$
and use a modification of $h$ to one with bounded second derivative as above.
 Then it follows that the vector field\refer{tilt} is Lipschitz on $\mathcal{X}$ and so defines a Lipschitz flow $\phi_F$ on $\mathcal{X}$.

 In this subsection we list some properties of the tilted FK model with phase space $\mathcal{X}$.

 Let $\phi_F^t\xbf$ denote a solution to\refer{tilt} with initial value $\xbf$ at $t=0$.
We remark that $\phi_F^t$ commutes with $T_{qp}$:
 \begeq\label{commu}
 T_{qp}\phi_F^t\xbf=\phi_F^t T_{qp}\xbf,\ \ \mbox{ for all }
 q,p\in\Z, \ \ \xbf\in\mathcal{X}, \mbox{ and }\ t\in\R.
 \eeq


 System\refer{tilt}  is said to be {\em monotone} if
 $\xbf\leq \ybf\Rightarrow \phi_F^t\xbf\leq \phi_F^t\ybf$ for  $t>0$
 for which both are defined. It is said to be {\em strictly monotone} if
 $\xbf < \ybf\Rightarrow \phi_F^t\xbf\ll \phi_F^t\ybf$ for  $t>0$
 for which both are defined. The condition $h_{12}(x,x')\leq -c<0$ guarantees
that system\refer{tilt} is strictly monotone~\cite{BM}.

In particular, by combining the results of the previous two paragraphs,
$$
\mathcal{X}_{MN} = \{\xbf \in \mathcal{X}\,|\, M\le x_{n+1}-x_n\le N,\ \forall n\in\Z \}
$$
 is invariant under $\phi_F^t$ for all $t>0$.  Thus, for all initial conditions in $\mathcal{X}_{MN}$ the modification of $h$ to achieve bounded second derivative has no effect on their trajectories.

Recall from \cite{BM} that the {\em width} $w(\xbf)$ of a configuration $\xbf \in \R^\Z$ is
$$
w(\xbf)=\sup \{ n^+(m)-n^-(m)\,|\, m \in \Z \},$$
where $n^-(m) = \max \{n \in \Z \,|\, T_{mn}\xbf \le \xbf \}$ and $n^+(m) = \min \{n\in \Z \,|\, T_{mn} \xbf \ge \xbf \}$. A configuration is Birkhoff if and only if its width equals $1$. If the width of $\xbf$ is bounded, then $w(\phi_F^t\xbf)\leq w(\xbf)$ for all $t>0$.

  Let
  $$
  \mathcal{X}_{\omega}=\{\xbf=(x_n)\in\mathcal{X}\,|\,\sup_{n\in\Z}
  |x_n-n\omega|<+\infty\}\mbox{ and }
  \mathcal{X}_\omega^1=\left\{
  \xbf\in\mathcal{X}_\omega\,|\,w(\xbf)\leq 1\right\}.
  $$
Assume $\xbf\in\mathcal{X}_\omega^1$. Then $\phi_F^t\xbf\in\mathcal{X}_\omega^1$ for all $t>0$. Hence
   $\mathcal{X}_\omega^1$ is forward
   invariant for\refer{tilt}.

   One of the most important quantities characterizing the  tilted FK model
  is the average velocity for each solution $\xbf(t)=(x_n(t))$
 of\refer{tilt} defined as the limit
 \begeq\label{nubar}
 \bar\nu=\lim_{{m-k\to\infty}\atop {t\to +\infty}}\frac{1}{t(m-k)}\int_0^t
 \sum_{n=k}^{m-1}\dot{x}_n(t)\mbox{d}t
 \eeq
 if it exists. It is a quantity separating sliding states and pinned states in particle chains~\cite{FM}.

It is proved in~\cite{QW1} the average velocity $\bar\nu$ exists in $\mathcal{X}_\omega^1$. It is independent of initial values but may depend on the tilt $F$ and the mean spacing $\omega$. So we denote it by $\nubar (\omega, F)$. Moreover, $\nubar(\omega,F)$ is nondecreasing in $F$.

  Define
  \begeq\label{dep_force}
  F_d(\omega)=\sup\{F\geq 0\,|\,\nubar(\omega,F)=0\}.
  \eeq
  Then $\nubar(\omega,F)=0$ for $0\leq F\leq F_d(\omega)$ and $\nubar(\omega,F)>0$ if $F>F_d(\omega)$.  The critical value $F_d(\omega)$ is called the depinning force~\cite{BM,BK,FM,Mac,PA}. The following proposition was proved in~\cite{BM}.
  \begpro\label{finite}
  For $0\leq F\leq F_d(p/q)$, there exist for\refer{tilt} equilibria of type $(p,q)$ and there are no equilibria of type $(p,q)$ while $F>F_d(p/q)$. Moreover, for $F>F_d(p/q)$ there is a time-periodic solution $\xbf(t)=(x_n(t))_{n\in\Z}$ of the gradient flow in the sense that $x_n(t+T) = x_n(t)+1$ for some $T>0$.  
  \epro

  Similar results for irrational $\omega$ were proved in~\cite{QW1}. Moreover, it was also proved (see also~\cite{WMWQ}) that $F_d(\omega)$ is continuous at irrational points, and it can be used as a criterion for the existence of invariant circles for the monotone twist map\refer{eq:apmap} with $F=0$: for irrational $\omega$, the area-preserving twist map $\Phi_0$  has an invariant circle with rotation number $\omega$ if and only if $F_d(\omega)=0$; for rational $\omega=p/q$ in lowest terms, $F_d(\omega)=0$ if and only if
there exists  an invariant circle on which the twist map is topologically conjugate to the rational rotation with rotation number $p/q$.


\subsection{Classic Aubry-Mather theory}\label{Sec:AM}

\begin{figure}[htbp] 
   \centering

\begin{tikzpicture}
  \draw [->](-3,0)--(3.8,0) node[anchor=north] {$n$};
  \foreach \x/\xtext in {-1/1,0/2,1/3,2/4,3/5}
  \draw (\x,0)--(\x,-2pt) node[below] {\scriptsize $\xtext$};
  \draw [->] (-2,0)--(-2,3.5) node[anchor=east] {$x$};
  \foreach \y/\ytext in {0.5/1,1/2,1.5/3,2/4,2.5/5,3/6}
  \draw (-2,\y)--(-2.05,\y) node[anchor=east] {\scriptsize $\ytext$};
  \draw (-2,0)--(3,2.5)node[anchor=west]{$\xbf$};
  \draw [blue](-2,0.5)--(3,3)node[above]{$T_{01}\xbf$};
  \draw [red] (-2,0.05)--(3,2.9)node[anchor=west]{$\ybf$};
  \foreach \x in {-2,-1,0,1,2,3}
  {\filldraw (\x,0.5*\x+1) circle (1pt);
  \filldraw [blue](\x,0.5*\x+1.5) circle (1pt);
  \filldraw [red] (\x,0.57*\x+1.18) circle (1pt);
  }

\end{tikzpicture}
   \caption{An example of Aubry diagram, showing a periodic state (in black) $\xbf$ of type $(1,1)$, its translate (blue) by $T_{01}$, and an advancing discommensuration (red) $\ybf$.}
   \label{fig:Aubrydgm}
\end{figure}
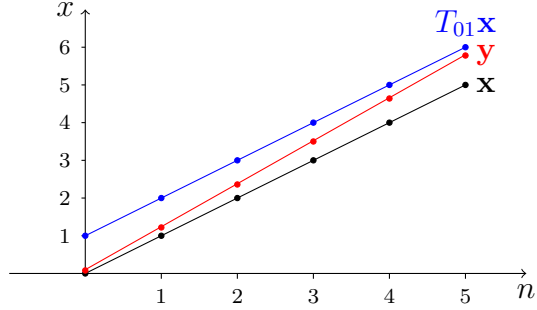

It is convenient to plot configurations $\xbf=(x_n)_{n \in \Z}$ in the plane of $(n,x_n)$, called the {\em Aubry diagram}, and to join their points by straight lines.  See Figure~\ref{fig:Aubrydgm}.
Then $T_{qp}$ is simply a translation by the vector $(q,p)$.

\begin{figure}[htbp] 
   \centering
   \begin{tikzpicture}[scale=0.8]
   \foreach \x in {0,1,...,11}
   \foreach \y in {0,-2.5}
   {
   \draw (-5.75+\x,1.25+\y)..controls(-5.6+\x,2.5+\y) and (-5.4+\x,2.5+\y)..(-5.25+\x,1.25+\y);
   \draw (-5.25+\x,1.25+\y)..controls(-5.1+\x,0+\y) and (-4.9+\x,0+\y)..(-4.75+\x,1.25+\y);
   }
   \filldraw [red] (-4.95,0.5) circle(1.5pt);
   \filldraw [red] (-3.87,0.8) circle(1.5pt);
   \filldraw [red] (-2.81,1.1) circle(1.5pt);
   \filldraw [red] (-1.8,1.4) circle(1.5pt);
   \filldraw [red] (-0.74,1.8) circle(1.5pt);

   \filldraw [red] (4.95,0.5) circle(1.5pt);
   \filldraw [red] (3.87,0.8) circle(1.5pt);
   \filldraw [red] (2.81,1.1) circle(1.5pt);
   \filldraw [red] (1.8,1.4) circle(1.5pt);
   \filldraw [red] (0.74,1.8) circle(1.5pt);

   \filldraw [blue] (-4.95,-2) circle(1.5pt);
   \filldraw [blue] (-4.05,-2) circle(1.5pt);
   \filldraw [blue] (-3.13,-1.7) circle(1.5pt);
   \filldraw [blue] (-2.19,-1.4) circle(1.5pt);
   \filldraw [blue] (-1.2,-1.1) circle(1.5pt);
   \filldraw [blue] (-0.26,-0.7) circle(1.5pt);

   \filldraw [blue] (4.95,-2) circle(1.5pt);
   \filldraw [blue] (4.05,-2) circle(1.5pt);
   \filldraw [blue] (3.13,-1.7) circle(1.5pt);
   \filldraw [blue] (2.19,-1.4) circle(1.5pt);
   \filldraw [blue] (1.2,-1.1) circle(1.5pt);
   \filldraw [blue] (0.26,-0.7) circle(1.5pt);

   \draw [red] (-4.95,0.5)--(-3.87,0.8)--(-2.81,1.1)--(-1.8,1.4)--(-0.74,1.8)--
   (0.74,1.8)--(1.8,1.4)--(2.81,1.1)--(3.87,0.8)--(4.95,0.5);

   \draw [blue] (-4.95,-2)--(-4.05,-2)--(-3.13,-1.7)--(-2.19,-1.4)--(-1.2,-1.1)--
   (-0.26,-0.7)--(0.26,-0.7)--(1.2,-1.1)--(2.19,-1.4)--(3.13,-1.7)--(4.05,-2)--(4.95,-2);

    \end{tikzpicture}


   \caption{Sketch of advancing (above) and retreating (below) discommensurations of type $(1,1)$ in the potential of a standard FK model.}
   \label{fig:discomms}
\end{figure}
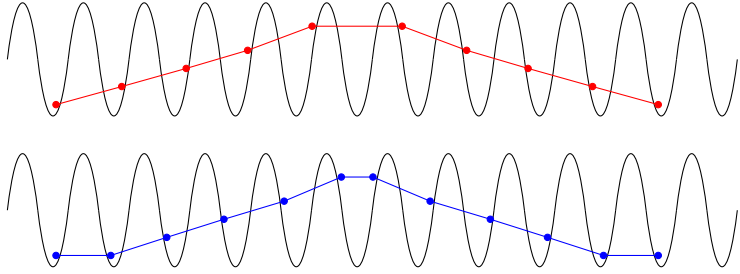

Among the equilibria of a non-tilted FK chain, particular importance can be given to the minimisers.  These are the configurations $\xbf$ such that for all $m<n-1$, the sum $W_{mn} = \sum_{i=m}^{n-1} h(x_i,x_{i+1})$ is minimum subject to $x_m,x_n$ fixed.
The theory of minimisers for non-tilted FK chains was established by Aubry \& Le Daeron \cite{ALD} and Mather \cite{Mat1}.  Surveys appear in \cite{Ban}, Ch.13 of \cite{KH}, \cite{MS} and \cite{Ba}.
The main results are that:
\begin{itemize}
\item every minimiser is Birkhoff. 
\item for each $p\in\Z, q\in\N$, there is a minimiser with $T_{qp}\xbf = \xbf$, called {\em type} $(p,q)$.
\item for $k\ge 2$, every minimiser of type $(kp,kq)$ is a minimiser of type $(p,q)$, so without loss of generality we take $p, q$ to be coprime.
\item the set $M_{p/q}$ of minimisers of type $(p,q)$ is totally ordered and closed.
\item if $M_{p/q}$ is not homeomorphic to $\R$ then for each gap in $M_{p/q}$, with endpoints denoted $\xbf^- < \xbf^+$, there is a minimiser $\xbf$ with $|x_n-x^-_n| \to 0$ as $n \to -\infty$ and $|x_n-x^+_n| \to 0$ as $n \to +\infty$, called an advancing discommensuration; similarly there is a minimiser with $|x_n-x^+_n| \to 0$ as $n \to -\infty$ and $|x_n-x^-_n| \to 0$ as $n \to +\infty$, called a retreating discommensuration.  See Figure~\ref{fig:discomms}.
\item denoting the sets of advancing and retreating discommensurations of mean spacing $p/q$ by $M_{p/q\pm}$ respectively, the sets $M_{p/q+}\cup M_{p/q}$ and $M_{p/q-}\cup M_{p/q}$ are totally ordered and closed in product topology (pointwise convergence).
\item for each irrational $\omega \in \R$ there is a minimiser with mean spacing $\omega$.
\item the set $M_\omega$ of minimisers with mean spacing $\omega$ irrational is totally ordered and closed in product topology.
\item the subset $M_\omega^r$ of recurrent minimisers is closed in product topology and its quotient by $T_{01}$ is homeomorphic to a circle or a Cantor set.
\end{itemize}

For tilted FK chains, there are no minimisers of bounded spacing $x_{n+1}-x_n$.  To see this, given $s >0$, let
$$E_1(s) = \sup | h(x,x'+1)-h(x,x')|, \ \ E_2(s) = \sup |h(x+1,x')-h(x,x')|$$ over $|x'-x| \le s$.  These are finite by periodicity and continuity of $h$.  Then the energy of any long enough segment $x_0,\ldots, x_L$ can be reduced by moving it to $x_0, x_1+1,\ldots, x_{L-1}+1,x_L$.  The change in energy is at most $E_1(s)+E_2(s)-(L-2)F$, which is negative if $L$ is large enough.
There may be Birkhoff equilibria, however.  Indeed, non-degenerate minimisers for $F=0$ persist as Birkhoff equilibria for $F$ small.

Given $s >0$, when $F$ is large enough there are no equilibria of spacing less than or equal to $s$.  To see this, let
$$K_1(s) = \sup h_2(x,x'), \ \ K_2(s) = \sup h_1(x,x')$$
over $|x'-x| \le s$.  Equilibrium requires $h_2(x_{n-1},x_n)+h_1(x_n,x_{n+1}) = F$.  So there are no equilibria with spacing bounded by $s$ if $F > K_1(s)+K_2(s)$.

\newsection{Invariant ordered circles of type $p/q\pm$}\label{sec:p/q+}
In this section we recall some definitions and results on invariant ordered circles (IOCs) and construct IOCs of type $p/q\pm$. Let $\one$ denote the element in $\mathcal{X}$ with each component being $1$.

\begde\label{defOC}
  The image $\ell$ of a continuous function $\tau$:  $\R\mapsto\mathcal{X}$ is called an {\em ordered circle} if $\tau(s_1)< \tau(s_2)$ for $s_1<s_2$ and $\tau(s+1)=\tau(s)+\one$. It is said to be a {\em strictly ordered circle} if
   furthermore $\tau(s_1)\ll \tau(s_2)$ for $s_1<s_2$. The function $\tau$ is called a parametrization of $\ell$.
  \ede

So an {ordered circle} is an ordered closed subset of $\R^\Z$ in product topology, homeomorphic to $\R$, and invariant under $T_{01}$. Its quotient by $T_{01}$ is homeomorphic to a circle.

\begde\label{defIOC}
We say that an ordered circle $\ell$ is a
 {\em translation-invariant ordered circle} (TOC) if $T_{qp}\ell=\ell$ for all $q$, $p\in\Z$.
An {\em invariant ordered circle} (IOC) is a TOC which is in addition invariant under the gradient flow $\phi_F^t$ of\refer{tilt} for all $t\in\R$.
\ede

Note that  each element in  a TOC $\ell$ is Birkhoff and has a mean spacing $\omega$, and hence $\ell\subset \mathcal{X}_{\omega}^1$.

For each $q \in \N, p \in \Z$, \cite{QW1} proved existence of an IOC $\ell_{p/q}$ of type $(p,q)$, i.e.~such that $T_{qp}\xbf = \xbf$ for all $\xbf \in \ell_{p/q}$.  By taking a limit point as $p/q \to \omega$ irrational, they proved existence of an IOC of mean spacing $\omega$ (their proof required addition of the assumption that the second derivative of $h$ is bounded, but as we explained at the beginning of Sec.2, this can be achieved by a simple modification).

\begpro\label{exisSOC} \cite{QW1}
  For each $\omega\in\R$ and each $F\geq 0$, there exists a strictly ordered circle $\ell\subset\mathcal{X}_{\omega}^1$ invariant for $\{T_{qp}\}_{q,p \in \Z}$ and $\phi_F^t$ for all $t\in\R$.
  \epro

\noindent{\bf Remark:}\
For each IOC $\ell\in\mathcal{X}_\omega^1$, there exists a lift of a circle homeomorphism $G$ with rotation number $\omega$ such that
$$
\ell=\{\xbf=(x_n)\,|\,x_n=G^n(s),\ n\in\Z,\ s\in\R\},
$$
owing to the invariance of $\ell$ under $T_{qp}$ for all $q, p\in\Z$, see~\cite{QW1}.
\vskip 1ex

It is useful to recall the proof of Proposition 3.3.

\vskip 1ex
\noindent{\bf Outline of proof of Prop.~\ref{exisSOC}}:
First \cite{QW1} consider the existence of IOC in
$$
\mathcal{X}_{p,q}^1=\{\xbf\in \mathcal{X}_{p/q}^1\,|\,T_{qp}\xbf=\xbf\},
$$
which is a finite-dimensional space. The Schauder fixed point theorem is used to prove the existence of an IOC $\ell_{p/q}$ in $\mathcal{X}_{p,q}^1$ for\refer{tilt}. Then
the $\ell_{p/q}$ modulo $T_{01}$ are compact subsets of $\mathcal{X}/T_{01}$.  For $p/q \to \omega$, they lie in a bounded subset of $\mathcal{X}/T_{01}$.  Thus in Hausdorff topology on compact subsets, they have a limit point.  Recall that two compact subsets are considered close in Hausdorff topology if each lies in a small neighbourhood of the other.  Lift such a limit point and call it $\ell_\omega$.

The states of $\ell_\omega$ are limit points in product topology of states with mean spacing $p/q$ and width at most 1, so they have mean spacing $\omega$ and width at most 1.

A key point is that since the $\ell_{p/q}$ are invariant under the gradient flow on $\mathcal{X}$, any limit point $\ell_\omega$ is also invariant under the gradient flow.

As a limit of ordered circles, $\ell_\omega$ is  ordered.  Using invariance under the gradient flow, we deduce it is strictly ordered, as follows.    Take two distinct states $\xbf, \ybf \in \ell_\omega$.  Flow backwards for some time $t$.  By invariance we obtain two distinct states of $\ell_\omega$.  They are weakly ordered.  Flow forwards by $t$.  By strict monotonicity of the gradient flow we obtain that $\xbf,\ybf$ are strictly ordered.

The set $\ell_\omega$ is also invariant under all the translations.  So we deduce it is an IOC of mean spacing $\omega$. \qed

We remark that if $0\leq F\leq F_d(\omega)$, then there is at least one equilibrium of\refer{tilt} on $\ell$ obtained in the above proposition. If $F>F_d(\omega)$, then the IOC $\ell$ obtained in the above proposition is a periodically sliding solution $x_n(t)=X(n\omega+vt+\alpha)$, where $X$ is called ``dynamical hull function'', and $v>0$. Therefore, for each $\xbf=(x_n)\in\ell$, $\dot{x}_n>0$. 

\begde\label{def-p/q}
An IOC $\ell$ is said to be of type $p/q\pm$ if each $\xbf \in \ell$ has mean spacing $p/q$, $T_{qp}\xbf \le \xbf$, respectively, $\ge \xbf$, and there exists $\xbf \in \ell$ for which $T_{qp}\xbf \ne \xbf$.
\ede

To construct IOCs of types $p/q\pm$, it
is enough to treat the case $p/q+$.  For the case $p/q-$, we can consider $\tilde{h}(x,x') = h(x',x)$ instead of $h$, which changes type $p/q-$ to type $(-p)/q+$.  Furthermore, one can add $F(x'-x)$ to $\tilde{h}$, which does not change the gradient flow.  It puts $h_F$ back into the form\refer{eq:hF}.

\begth\label{construction-IOC-p/q}
Let $0< \bar{F}<F_d(p/q)$. Then there exists an IOC of type $p/q+$ for\refer{tilt} with $F=\bar{F}$.
\eeth
\pf
Choose $\omega_n\searrow p/q$, $F_n\to \bar{F}$, as $n\to\infty$, and $\ell_n$ be an IOC with mean spacing $\omega_n$ for\refer{tilt} with $F=F_n$.
To construct IOCs of type $p/q+$, we follow a strategy similar to~\cite{QW2}. The quotients of
$\ell_n$ by $T_{01}$ have a limit point in Hausdorff topology on compact subsets of $\mathcal{X}/T_{01}$. Lift one such limit point and call it $\ell_{p/q+}$.
Each $\xbf \in \ell_{p/q+}$ has mean spacing $p/q$ and width at most 1, because it is a limit in product topology of states with mean spacing $\omega_n$ of width at most 1.
Since $T_{qp}\ybf<\ybf$ for each $\ybf\in\ell_n$ due to $\omega_n>p/q$, then for each $\xbf\in\ell_{p/q+}$ which is a limit in product topology of $\ybf^n\in\ell_n$, it holds that $T_{qp}\xbf\leq\xbf$.
It remains to show that there exists $\xbf\in\ell_{p/q+}$ such that $T_{qp}\xbf\ne \xbf$.
Note that we can choose $\omega_n$ such that $\lim_{n\to\infty}F_d(\omega_n)=\limsup_{\omega\searrow p/q}F_d(\omega)$.

\vspace{2mm}

Case (i)\ $0< \bar{F}\leq \limsup_{\omega\searrow p/q}F_d(\omega)$. If the righthand inequality is strict, take $F_n=\bar{F}$. If it is an equality, let $F_n=F_d(\omega_n)$. In either case, we have $F_n\leq F_d(\omega_n)$ for $n$ large enough. It then follows that
 there exists at least one equilibrium on $\ell_n$ for\refer{tilt} with
$F=F_n$, implying the existence of at least one equilibrium on $\ell_{p/q+}$ for\refer{tilt} with $F=\bar{F}$.

Moreover, the set of equilibria on $\ell_{p/q+}$ contains at least one of type $(p,q)$ because we can take a limit point in product topology of translates of any equilibrium on it by $T_{qp}^k$,
$k\to\infty$. The set $M$ of equilibria of type $(p,q)$ on $\ell_{p/q+}$ is closed. The IOC $\ell$ can not consist entirely of equilibria,
because that would give a rotational\footnote{homotopically non-trivial} invariant circle for the corresponding area-preserving map $\Phi_{\bar{F}}$ determined in\refer{eq:apmap},
but the map has net flux $\bar{F} \ne 0$ so it has no rotational invariant circles.

We now show that in each gap of $M$ on $\ell_{p/q+}$ there is a $p/q+$ equilibrium, so in particular there exists $\xbf$ in the gap such that $T_{qp}\xbf<\xbf$.
This is a consequence of the same argument as was used in~\cite{ALD} for the untilted case; see also the approach in Section 5 in~\cite{Ban} that we adopt here. Let
$\xbf^-<\xbf^+$ denote the endpoints of the gap and
$\xbf^n\in\ell_n$ be an equilibrium of\refer{tilt} with $F=F_n$. Let $\eps=\min_{i\in\Z}(x^+_i-x^-_i)>0$. Then we can choose appropriate translates of the $\xbf^n$, denoted by $\tilde{\xbf}^n$, such that for some $j\in\{1,2,\cdots,q\}$,
$$
\tilde{x}_{j-1}^n\leq x_{j-1}^-+\eps/2, \ \mbox{ and }
\tilde{x}_j^n>x_j^-+\eps/2.
$$
There is a convergent subsequence of $\tilde{\xbf}^n$.  Its limit $\xbf\in \ell_{p/q+}$ because $\ell_n$ converge to $\ell_{p/q+}$.  $\xbf$ is an equilibrium of\refer{tilt} with $F=\bar{F}$.  Because it has points in the gap between $\xbf^-$ and $\xbf^+$ and all three belong to the same IOC, it satisfies $\xbf^-<\xbf<\xbf^+$, $T_{qp}\xbf<\xbf$, and $T_{qp}^k\xbf\to \xbf^-$ as $k\to +\infty$, $T_{qp}^k\xbf\to \xbf^+$ as $k\to -\infty$. 

\vspace{2mm}

Case (ii)\ $\limsup_{\omega\searrow p/q}F_d(\omega)\leq \bar{F}<F_d(p/q)$.  If the left inequality is strict, take $F_n=\bar{F}$. If it is an equality, let $F_n=F_d(\omega_n)$. In either case, we have $F_n\geq F_d(\omega_n)$ for $n$ large enough. It then follows that for each $\ybf\in\ell_n$, $\dot{\ybf}\geq\zero$ for\refer{tilt} with $F=F_n$ (see~\cite{QW1}), and hence for each $\xbf\in\ell_{p/q+}$, which is a limit in product topology of $\ybf^n\in\ell_n$, $\dot{\xbf}\geq\zero$ for\refer{tilt} with $F=\bar{F}$.

Let $\hat{F} \in (\bar{F}, F_d(p/q))$. Then there exists an equilibrium $\hat{\xbf}$ of type $(p,q)$ for\refer{tilt} with $F=\hat{F}$ thanks to Proposition~\ref{finite} (not necessarily on $\ell_{p/q+}$). Assume all the configurations on $\ell_{p/q+}$ are of type $(p,q)$, i.e., $T_{qp}\xbf=\xbf$ for $\xbf\in\ell_{p/q+}$. Let
$$
\wbf=\sup\{\xbf\in\ell_{p/q+}\,|\,\xbf\leq\hat{\xbf}\}.
$$
Then $\wbf\in\ell_{p/q+}$ and is of type $(p,q)$. Next, there exists $n\in\Z$ such that
$w_{n-1}\leq \hat{x}_{n-1}$, $w_n=\hat{x}_n$, and $w_{n+1}\leq \hat{x}_{n+1}$ due to the continuity of $\ell_{p/q+}$.
It then follows from  the monotonicity that
$$
\begin{array}{lll}
-h_2(\hat{x}_{n-1},\hat{x}_n)-h_1(\hat{x}_n,\hat{x}_{n+1})+\hat{F}& > & -h_2(\hat{x}_{n-1},\hat{x}_n)-h_1(\hat{x}_n,\hat{x}_{n+1})+\bar{F}\\[4pt]
 & \geq & -h_2(w_{n-1},w_n)-h_1(w_n,w_{n+1})+\bar{F}\\[4pt]
  &= &\dot{w}_n\geq 0,
 \end{array}
$$
 a contradiction to the fact that $\hat{\xbf}$ is an equilibrium of the tilted FK chain\refer{tilt} with $F=\hat{F}$. Consequently, there exists $\xbf\in\ell_{p/q+}$ such that $T_{qp}\xbf\neq \xbf$.
\emark

\noindent{\bf Remark:}\ For $0< \bar{F}\leq \limsup_{\omega\searrow p/q}F_d(\omega)$, from the above proof of case (i) we have actually obtained that in each gap of $M$, there is a $p/q+$ equilibrium, where $M$ denotes the set of equilibria on $\ell_{p/q+}$ of type $(p,q)$. If $\limsup_{\omega\searrow p/q}F_d(\omega_n)\leq \bar{F}<F_d(p/q)$, then there exists an IOC $\ell_{p/q+}$ such that $\dot{\xbf}\geq\zero$ for each $\xbf\in\ell_{p/q+}$, and in at least one gap of $M$, there is a configuration of type $p/q+$, which is not necessarily an equilibrium.

\vspace{2mm}

Let
$$
\mcalX_{p,q}=\{\xbf\in\R^\Z\,|\, T_{qp}\xbf=\xbf\},
$$
and
$$
W_{p,q}(\xbf)=\sum_{i=0}^{q-1}h_F(x_i,x_{i+1}),\ \ \mbox{ for }\ \xbf=(x_i)\in\mcalX_{p,q}.
$$
Then $\xbf\in\mcalX_{p,q}$ is an equilibrium of\refer{tilt} if and only if it is a critical point of $W_{p,q}$.

Let $\xbf\in\mcalX_{p,q}$ be a  critical point of $W_{p,q}$. By the Morse index of $W_{p,q}$ at $\xbf$ we mean the number of negative eigenvalues of the Hessian matrix of $W_{p,q}$ at $\xbf$.

Let $\ell_F$ be an IOC of type $p/q+$ obtained by Theorem~\ref{construction-IOC-p/q} for\refer{tilt} with $F\in (0,F_d(p/q))$. In what follows, we shall show that on $\ell_F$ there are no  equilibria of type $(p,q)$ for\refer{tilt} with positive Morse index, which implies that if the equilibria of type $(p,q)$ on $\ell_F$ are non-degenerate, then they are local minima for $W_{p,q}$. We use results of~\cite{An1}.

We temporarily denote $\limsup_{\omega\searrow p/q}F_d(\omega)$ by $F_d(p/q+)$, the justification of which will be given by the proof of Theorem A.

\begth\label{positive-index}
Let $F\in (0, F_d(p/q))$ and $\ell_F$ be an IOC  of type $p/q+$ obtained by Theorem~\ref{construction-IOC-p/q}. In addition, assume $\dot{\zbf}\geq\zero$ for each $\zbf\in\ell_F$ if $F_d(p/q+)<F<F_d(p/q)$. Then on $\ell_F$ there are no  equilibria for\refer{tilt} of type $(p,q)$ with positive Morse index.
\eeth
\noindent \pf
Assume to the contrary there exists an equilibrium $\xbf^*=(x_n^*)\in\ell_F$ of type $(p,q)$ with positive Morse index.

We remark that from equilibria $\xbf=(x_n)$ of\refer{tilt} we can define another planar map $\tilde{\Phi}_F$: $(x_n,x_{n+1})\mapsto (x_{n+1},x_{n+2})$ which is conjugate to $\Phi_F$ defined by\refer{eq:apmap}. In this case $(x_0^*,x_1^*)$ is a $(p,q)$-periodic point of $\tilde{\Phi}_F$. Denote the Jacobian matrix $D\tilde{\Phi}_F^q$ at $(x_0^*,x_1^*)$ by $A$ for simplicity (note that $A$ is identical to $M(0)$ of Section 3 in~\cite{An1}). The following conclusion is derived by the discussion of Section 3 in~\cite{An1} (see also~\cite{MM}).

Since $\xbf^*$ has positive Morse index, then there is a positive integer $N\in\N$, depending on the Jacobian matrix $A$, such that for sufficiently small $\eps_0>0$, if there exists an equilibrium $\xbf=(x_n)$ of\refer{tilt} with $|x_n-x_n^*|<\eps_0$ for $0\leq n\leq N$, then there exist $i,j\in\Z$ with $0<i<j<N$ such that
$$
(x_i-x_i^*)(x_j-x_j^*)<0.
$$

Let $0<F\leq F_d(p/q+)$. Then from the remark of Theorem~\ref{construction-IOC-p/q}, there exists a $p/q+$ type equilibrium $\xbf=(x_n)\in\ell_F$ in a gap with the right endpoint being $\xbf^*$, implying $\xbf\ll\xbf^*$ (the proof is the same if $\xbf^*$ is degenerate and all the configurations in the gap are equilibria of type $(p,q)$). Meanwhile, $T_{qp}^k\xbf\to\xbf^*$ as $k\to -\infty$ implies $|\tilde{x}_n-x_n^*|<\eps_0$ for $0\leq n\leq N$, where $\tilde{\xbf}=(\tilde{x}_n)=T_{qp}^k\xbf$ for some $k\in\Z$. This is a contradiction since we have $\tilde{x}_n<x_n^*$ for all $n\in\Z$.

Let $F_d(p/q+)<F<F_d(p/q)$ and $\ell_F$ be an IOC of type $p/q+$ of\refer{tilt} satisfying $\dot{\zbf}\geq \zero$ for all $\zbf\in\ell_F$.

Let $\zbf=(z_n)\in\ell_F$ be in the gap with the right endpoint being $\xbf^*$, implying $\zbf\ll\xbf^*$. If all the configurations in the gap are equilibria of type $(p,q)$, or there is  an equilibrium of type $p/q+$ in the gap,  then we have a contradiction, as in the case $0<F\leq F_d(p/q+)$. So we assume $\dot{\zbf}\gg\zero$.

Let $x_0=z_0$ and $x_1=z_1$. Then solving\refer{=F} to obtain $x_2,x_3,\cdots$, and $x_{-1},x_{-2},\cdots$., we have an equilibrium $\xbf=(x_n)$ of\refer{tilt}. Note that by the condition $h_{12}(x,x')<-c$ we have $x_2<z_2$ and $x_{-1}<z_{-1}$. Indeed, since
$$
-h_2(z_0,z_1)-h_1(z_1,z_2)+F>0,\ \ \
-h_2(z_0,z_1)-h_1(z_1,x_2)+F=0,
$$
and $-h_1$ is strictly increasing with respect to its second variable, then $x_2<z_2$. The proof for $x_{-1}<z_{-1}$ is similar using the fact that $-h_2$ is strictly increasing with respect to its first variable.

Furthermore, we have the conclusion that $x_n\leq z_n$ for all $n\in\Z$, i.e., $\xbf\leq \zbf$. Indeed, assume to the contrary there exists some $j>2$ ($j<-1$ is similar) such that $x_j\geq z_j$. Without loss of generality we assume $j=3$. Note that the IOC $\ell_F$ has a parametrization $\tau(s)$. Let $\tau(s_2)=\zbf$. Then there exists $s_1<s_2$ such that $w_2=x_2<z_2$, where $\wbf=(w_n)=\tau(s_1)\in\ell_F$. Meanwhile we have $w_1<z_1=x_1$ and $w_3<z_3\leq x_3$ due to $\tau(s_1)\ll\tau(s_2)$. Since $\dot{\wbf}\geq \zero$, then
$$
-h_2(w_1,w_2)-h_1(w_2,w_3)+F\geq 0, \ \ \mbox{ hence }\
-h_2(x_1,x_2)-h_1(x_2,x_3)+F>0,
$$
a contradiction. Consequently, it follows $\xbf\leq \zbf\ll\xbf^*$.

On the other hand, we can choose $z_0$ and $z_1$ sufficiently close to $x_0^*$ and $x_1^*$ respectively, such that $|x_n-x_n^*|<\eps_0$ for $0\leq n\leq N$. Then we have $0<i<j<N$ with $(x_i-x_i^*)(x_j-x_j^*)<0$, a contradiction to $\xbf\ll\xbf^*$.

Therefore, there are no  equilibria of type $(p,q)$ with positive Morse index on $\ell_F$.\emark

We remark that the equilibria of type $p/q+$ have the similar property if they lie on an IOC of type $p/q+$ for $F_d(p/q+)\leq F<F_d(p/q)$. More precisely,
let $\ell_F$ be an IOC of type $p/q+$ of\refer{tilt} with the property $\dot{\zbf}\geq \zero$ for each $\zbf\in\ell_F$, $\xbf^*=(x_n^*)\in\ell_F$ be an equilibrium of type $p/q+$, and
$$
W_{0,n+1}^*(y)=h_F(x_0^*,y_1)+\sum_{i=1}^{n-1}h_F(y_i,y_{i+1})+h_F(y_n,x_{n+1}^*),
$$
where $y=(y_1, \cdots,y_n)\in\R^n$. Consider the truncated system
$$
\dot{y}=-\nabla W_{0,n+1}^*(y),
$$
i.e.,
\begeq\label{truncate}
\left\{
\begin{array}{lll}
\dot{y}_1 &= & -h_2(x_0^*,y_1)-h_1(y_1,y_2)+F,\\[4pt]
\dot{y}_i &= & -h_2(y_{i-1},y_i)-h_1(y_i,y_{i+1})+F,\ \ i=2,\cdots,n-1,\\[4pt]
\dot{y}_n &= & -h_2(y_{n-1},y_n)-h_1(y_n,x_{n+1}^*)+F.
\end{array}
\right.
\eeq
Note that\refer{truncate} is strictly monotone on $\R^n$ and $y^*=(x_1^*,\cdots,x_n^*)$ is an equilibrium of\refer{truncate}.

We say that the equilibrium $\xbf^*$ of\refer{tilt} has positive Morse index from $0$ to $n+1$ if $y^*$ has positive Morse index as an equilibrium of\refer{truncate}, i.e., the $n\times n$ matrix $-D^2W_{0,n+1}^*$ has at least one positive eigenvalue at $y^*$. Similarly for $l,m\in\Z$ with $l<m-1$, we can give the definition that $\xbf^*$ has positive Morse index from $l$ to $m$ by considering the truncated energy function $W_{l,m}^*(y)$, in which $y\in\R^{m-l-1}$.

\begth\label{posi-index-p/q+}
Let $F_d(p/q+)\leq F<F_d(p/q)$ and $\ell_F$ be an IOC of type $p/q+$ with the property $\dot{\zbf}\geq\zero$ for each $\zbf\in\ell_F$. Then there are no equilibria of type $p/q+$ on $\ell_F$ with positive Morse index from some $l$ to $m$ ($l<m-1$).
\eeth
\noindent\pf
Without  loss of generality, we assume to the contrary $\xbf^*=(x_n^*)\in\ell_F$ is an equilibrium of type $p/q+$ with positive Morse index from $0$ to $n+1$, i.e., $y^*$ has positive Morse index as an equilibrium of the truncated system\refer{truncate}. Then from the Perron-Frobenius Theorem it follows that the matrix $-D^2W_{0,n+1}^*(y^*)$ has a positive eigenvalue corresponding to a positive eigenvector $\xi=(\xi_i)\in\R^n$ with $\xi_i>0$ for $i=1,\cdots,n$. Then (see Section 5, Chapter 2 in~\cite{Sm}) there is a neighbourhood $U$ of $y^*$ such that if $y(0)\in U$ and $y(0)\ll y^*$, then $y(t)\ll y(0)$ for all $t>0$, where $y(t)$ is the solution of\refer{truncate} with initial value $y(0)$ at $t=0$.

On the other hand, we can choose $\zbf=(z_n)\in\ell_F$, $\zbf\ll\xbf^*$, sufficiently close to $\xbf^*$, implying $z_i$ close to $x_i^*$ for $i=1,\cdots,n$. Note that $\dot{z}_i\geq 0$, $i\in\Z$ for\refer{tilt}. Comparing\refer{tilt} with\refer{truncate}, we deduce that $\dot{z}_i\geq 0$, $i=1,\cdots, n$ for\refer{truncate} owing to the facts $z_0\leq x_0^*$ and $z_{n+1}\leq x_{n+1}^*$. If we set $y(0)=(z_1,\cdots,z_n)$, then the solution $y(t)$ of\refer{truncate} with initial value $y(0)$ at $t=0$ has the property $y(t)\geq y(0)$ for $t>0$ since\refer{truncate} is strictly monotone. This is a contradiction to $y(t)\ll y(0)$ for $t>0$ obtained by conclusions in~\cite{Sm}.\emark

We say that an equilibrium $\xbf^*$ of type $p/q+$ is degenerate if $0$ is in the spectrum of $A=-D^2W(\xbf^*)$ which is a bounded operator on the Hilbert space $\ell^2$.

\begth\label{degenerate-p/q+}
Let $F_d(p/q+)\leq F<F_d(p/q)$ and $\ell_F$ be an IOC of type $p/q+$ with the property $\dot{\zbf}\geq\zero$ for each $\zbf\in\ell_F$. Then all equilibria of type $p/q+$ on $\ell_F$ are degenerate.
\eeth
\noindent\pf
Let $\xbf^*\in\ell_F$ be an equilibrium of type $p/q+$ and $A=-D^2W(\xbf^*)$.  Then for $\eta=(\eta_i)\in\ell^2$,
\begeq\label{Aeta}
(A\eta)_i=\alpha_i\eta_{i-1}+\beta_i\eta_i+\alpha_{i+1}\eta_{i+1},
\eeq
where $\alpha_i=-h_{12}(x_{i-1}^*,x_i^*)$, $\beta_i=-h_{22}(x_{i-1}^*,x_i^*)-h_{11}(x_i^*,x_{i+1}^*)$. We shall show that each $\lambda>0$ is in the resolvent set of $A$.

Let $\xi=(\xi_i)\in\ell^2$. We show that the equation $(\lambda I-A)\eta=\xi$ has a unique solution in $\ell^2$, where $I$ denotes the identity of $\ell^2$ and $\lambda>0$.

For $n\geq 1$, let $W_{-n,n}^*(y)$ be defined as in the proof of Theorem~\ref{posi-index-p/q+} and $A_n=-D^2W_{-n,n}^*(y^*)$, where $y^*=(x_{-n+1}^*,\cdots,x_{n-1}^*)$. Then $A_n$ is a symmetric matrix of order $2n-1$. From Theorem~\ref{posi-index-p/q+} we know that $A_n$ has no positive eigenvalues.

Let $\bar{\xi}^n=(\xi_{-n+1},\cdots,\xi_0,\cdots,\xi_{n-1})\in\R^{2n-1}$. Then there exists $\bar{\eta}^n\in\R^{2n-1}$ such that
$$
(\lambda I_n-A_n)\bar{\eta}^n=\bar{\xi}^n,
$$
where $I_n$ denotes the identity of the Euclidean space $\R^{2n-1}$ with the Euclidean norm $|\cdot|_{2n-1}$.

Let $\mu_{2n-1}\leq\cdots\leq\mu_1\leq 0$ denote the eigenvalues of $A_n$. There exists an orthonormal matrix $P_n$ of order $2n-1$ such that $P_n^{-1}A_nP_n=\mbox{diag}(\mu_{2n-1},\cdots,\mu_1)$. Let
$$
\tilde{\eta}^n=P_n^{-1}\bar{\eta}^n\  \ \mbox{ and }\ \ \tilde{\xi}^n=P_n^{-1}\bar{\xi}^n.
$$
Then
$$
P_n^{-1}(\lambda I_n-A_n)P_n\tilde{\eta}^n=\tilde{\xi}^n,
$$
hence
$$
\mbox{diag}(\lambda-\mu_{2n-1},\cdots,\lambda-\mu_1)\tilde{\eta}^n=\tilde{\xi}^n,
$$
implying
\begeq\label{etabar}
|\bar{\eta}^n|_{2n-1}=|\tilde{\eta}^n|_{2n-1}\leq \frac{1}{\lambda}|\tilde{\xi}^n|_{2n-1}=
\frac{1}{\lambda}|\bar{\xi}^n|_{2n-1}.
\eeq

Let $\hat{\eta}^n=(\hat{\eta}_i^n)\in\ell^2$ with $\hat{\eta}_i^n=\bar{\eta}_i^n$ for $i=-n+1,\cdots,n-1$, and $\hat{\eta}_i^n=0$ otherwise. Then $\{\hat{\eta}^n\}$ has a limit point, denoted by $\eta=(\eta_i)$, in the product topology, i.e., for each $i\in\Z$, $\hat{\eta}_i^n\to\eta_i$ as $n\to\infty$.

Let $\|\cdot\|_2$ denote the norm of $\ell^2$. Note that for each $k\geq 1$,
$$
\sum_{i=-k}^k\eta_i^2=\lim_{n\to\infty}\sum_{i=-k}^k(\hat{\eta}_i^n)^2\leq\lim_{n\to\infty}
\|\hat{\eta}^n\|_2^2=\lim_{n\to\infty}|\bar{\eta}^n|_{2n-1}^2\leq\lim_{n\to\infty}
\frac{1}{\lambda^2}
|\bar{\xi}^n|_{2n-1}^2\leq\frac{1}{\lambda^2}\|\xi\|_2^2.
$$
It then follows that $\eta\in\ell^2$ and $\|\eta\|_2\leq (1/\lambda)\|\xi\|_2$. Meanwhile, it is easy to check that $\eta$ satisfies $(\lambda I-A)\eta=\xi$.

  It remains to show that the solution is unique. Assume to the contrary the solution is not unique, which is equivalent to that $(\lambda I-A)\eta=\zero$ has a non-zero solution $\eta=(\eta_i)$ with $\|\eta\|_2=1$. Let
  $$
  \bar{\eta}^n=(\eta_{-n+1},\cdots,\eta_{n-1})\in\R^{2n-1}.
  $$
  Then from\refer{Aeta} and the calculation of $A_n$ we have
  $$
  (\lambda I_n-A_n)\bar{\eta}^n=\bar{\xi}^n=(-\alpha_{-n+1}\eta_{-n},0,\cdots,0,
  -\alpha_n\eta_n)\in\R^{2n-1}.
  $$
  Then
  $$
  |\bar{\xi}^n|_{2n-1}^2\leq \alpha^2(|\eta_{-n}|^2+|\eta_n|^2), \ \ \mbox{where }
  \ \alpha=\sup_{i\in\Z}\{|\alpha_i|\}<\infty.
  $$
  From\refer{etabar} it follows that
  $$
  |\bar{\eta}^n|_{2n-1}^2\leq \frac{1}{\lambda^2}|\bar{\xi}^n|_{2n-1}^2\leq
  \frac{\alpha^2}{\lambda^2}\left( |\eta_{-n}|^2+|\eta_n|^2\right).
  $$
  Choose $n$ large enough such that
  $$
  |\bar{\eta}^n|_{2n-1}^2\geq 1/2 \ \ \mbox{ and }\ \ \frac{\alpha^2}{\lambda^2}\left( |\eta_{-n}|^2+|\eta_n|^2\right)<1/2,
  $$
  yielding a contradiction.

  As a consequence, the solution for $(\lambda I-A)\eta=\xi$ is unique, and hence $\lambda>0$ is in the resolvent set of $A$, i.e., $\sigma(A)\subset (-\infty,0]$, where $\sigma(A)$ denotes the spectrum of $A$.

In what follows we shall show that $0\in\sigma(A)$. Assume to the contrary $0\not\in\sigma(A)$. Then $\sigma(A)\subset (-\infty,-c]$ for some $c>0$, implying that if $\xbf-\xbf^*\in\ell^2$ and $\|\xbf-\xbf^*\|_2$ is small enough, then $\|\xbf(t)-\xbf^*\|_2\to 0$ as $t\to +\infty$, where $\xbf(t)$ is a solution of\refer{tilt} with $\xbf(0)=\xbf$.

However, there exists $\xbf\in\ell_F$, $\xbf^*\ll\xbf$, $\|\xbf-\xbf^*\|_2$ small enough, but $\dot{\xbf}\geq\zero$, yielding a contradiction. Indeed, we need to verify that $\zbf-\xbf^*\in\ell^2$ for $\zbf\in\ell_F$ and $\xbf^*\ll\zbf\ll T_{qp}^{-1}\xbf^*$. Note that
there are two equilibria $\xbf^-$ and $\xbf^+$ of type $(p,q)$ such that $\xbf^-\ll\xbf^*\ll\xbf^+$ and $T_{qp}^{-k}\xbf^*\to \xbf^+$, $T_{qp}^k\xbf^*\to\xbf^-$ as $k\to +\infty$ in the product topology. Meanwhile,
$$
\sum_{i\in\Z}z_i-x_i^*\leq\sum_{i\in\Z}(T_{qp}^{-1}\xbf^*)_i-x_i^*=
\lim_{n\to +\infty}\sum_{i=1}^q(T_{qp}^{-n}\xbf^*)_i-(T_{qp}^n\xbf^*)_i
=\sum_{i=1}^q x_i^+-x_i^-<\infty.
$$
As a consequence, $\zbf-\xbf^*\in\ell^2$ for each $\zbf\in\ell_F$ with $\xbf^*\ll\zbf\ll T_{qp}^{-1}\xbf^*$.

Therefore, we  choose  $\xbf\in\ell_F$ with $\xbf^*\ll\xbf\ll T_{qp}^{-1}\xbf^*$ such that $\|\xbf-\xbf^*\|_2$ is small enough. Since $\dot{\xbf}\geq\zero$, the solution $\xbf(t)$ of\refer{tilt} with initial value $\xbf(0)=\xbf$ has the property $\xbf(t)\geq\xbf\gg \xbf^*$ for $t>0$, a contradiction to $\|\xbf(t)-\xbf^*\|_2\to 0$ as $t\to +\infty$. Then we arrive at the conclusion that $0\in\sigma(A)$ and hence $\xbf^*$ is degenerate.  \emark

\newsection{Proofs of main results}\label{sec:pf}
\subsection{Preparation}
To proceed, we need some more preparatory results.
\begde\label{glue}
Let $\ybf=(y_n)$, $\zbf=(z_n)\in\R^\Z$, and $x_n=y_n$ for $n\leq n_0$ and $x_n=z_n$ for $n>n_0$. We say that $\xbf=(x_n)$ is a $\delta$-gluing of $\ybf$ and $\zbf$ if $|y_n-z_n|<\delta$ for $n=n_0,\,n_0+1$. Similarly, we can define a $\delta$-gluing of more than two configurations.
\ede
Note that by continuity of $h_1$ and $h_2$ in\refer{tilt}, if $\ybf$ and $\zbf$ are equilibria of\refer{tilt} with bounded spacing, then for each $\eps>0$, a $\delta$-gluing $\xbf$ of $\ybf$ and $\zbf$ satisfies
$|\dot{x}_n|<\eps$ for\refer{tilt} if $\delta$ is small enough.  Indeed, $|\dot{x}_n| < C \delta$ on $\mathcal{X}_{MN}$ for $C=\sup \{|h_{12}(x,x')| : M\le x'-x \le N\}$.

\begin{figure}
  \begin{center}
  \begin{picture}(320,120)

  \put(0,0){\line(1,0){320}}
  \put(0,0){\line(0,1){120}}
 \put(320,0){\line(0,1){120}}
  \put(160,0){\line(0,1){120}}
 \put(0,120){\line(1,0){320}}

   \put(0,0){\line(4,1){160}}
  \put(0,80){\line(4,1){160}}
   \put(20,5){\circle*{2}}
   \put(40,10){\circle*{2}}
    \put(60,19){\circle*{2}}
    \put(80,34){\circle*{2}}
   \put(100,65){\circle*{2}}
   \put(120,95){\circle*{2}}
   \put(140,115){\circle*{2}}

  \qbezier(100,65)(120,110)(160,118)
  \qbezier(100,65)(80,20)(40,12)
  \put(120,25){\makebox(0,0)[rt]{$\ybf$}}
  \put(120,85){\makebox(0,0)[rt]{$\zbf$}}
  \put(40,95){\makebox(0,0)[rb]{$T_{q'p'}\ybf$}}

   \put(160,0){\line(4,1){160}}
  \put(160,40){\line(4,1){160}}
  \put(160,80){\line(4,1){160}}

 \qbezier(290,97.5)(300,115)(320,118)
  \qbezier(290,97.5)(280,70)(260,68)
\qbezier(230,37.5)(240,60)(260,63)
 \qbezier(230,37.5)(220,15)(200,12)

\put(180,5){\circle*{2}}
   \put(200,10){\circle*{2}}
    \put(220,22){\circle*{2}}
    \put(240,52){\circle*{2}}
   \put(260,65){\circle*{2}}
   \put(280,78){\circle*{2}}
   \put(300,115){\circle*{2}}

\put(280,25){\makebox(0,0)[rt]{$\ybf$}}
\put(200,95){\makebox(0,0)[rb]{$T_{q'p'}\ybf$}}
\put(280,65){\makebox(0,0)[rt]{$\wbf$}}
\put(305,100){\makebox(0,0)[rt]{${\zbf}''$}}
\put(240,35){\makebox(0,0)[rt]{${\zbf}'$}}

  \end{picture}
  \end{center}
  \caption{(a) The construction of $\xbf$, by translating copies of segments of  $\ybf$ and $\zbf$; (b) Construction with an intermediate equilibrium $\wbf$, using ${\zbf}', {\zbf}''$ (see Sec.~4.2(a)).
  }\label{Fig:x}
  \end{figure}
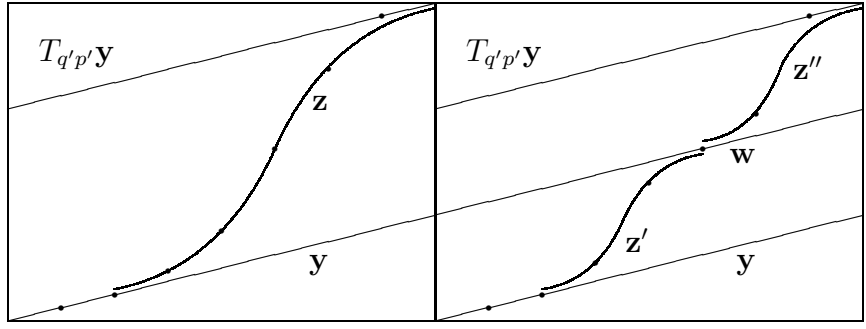

\begle\label{per-glue}
Let $\ybf=(y_j)$ be a $(p,q)$-periodic equilibrium of\refer{tilt} and $\ybf\ll T_{q'p'}\ybf$ for some $q'\in\N$ and $p'\in\Z$. Assume for each $\eps>0$, there exists a configuration $\zbf=(z_j)$ such that
$$
|\dot{z}_j|=|-h_2(z_{j-1},z_j)-h_1(z_j,z_{j+1})+F|<\eps,
$$
and
$$
|z_j-y_j|\to 0 \ \mbox{ as } j\to -\infty \ \mbox{ and } \ |z_j-(T_{q'p'}\ybf)_j|\to 0 \ \mbox{ as } j\to +\infty.
$$
Then there exists $N\in\N$ such that for each $n\geq N$, there is a periodic configuration
$\xbf=(x_j)$ of type $(np+p',nq+q')$ satisfying $|\dot{x}_j|<2\eps$ for all $j\in\Z$.
\ele
\noindent\pf
Since the translates $\{T_{qp}^m\zbf\,|\,m\in\Z\}$ of $\zbf$ also satisfy the assumptions of the lemma, we may  replace $\zbf$ by its translate $T_{qp}^m\zbf$ if necessary. Then
given $\delta>0$, there exists $N\in\N$ such that for $n\geq N$,
$$
|z_j-y_j|<\delta  \mbox{ for }\,j\leq q', \mbox{ and }
|z_j-(T_{q'p'}\ybf)_j|<\delta\,\mbox{ for }\,j\geq q'+nq-1.
$$
Let $P_n=np+p'$ and $Q_n=nq+q'$. We make a configuration $\xbf$ of type $(P_n,Q_n)$ as follows, see
Figure~\ref{Fig:x}(a).
Let
$$
\begin{array}{l}
x_j=y_j, \mbox{ for } j=0,1,\cdots,q'-1,\\[4pt]
 x_j=z_j,\mbox{ for }j=q',\cdots, q'+nq-1,\mbox{ and }\\[4pt]
x_j=(T_{q'p'}\ybf)_j, \mbox{ for }j=q'+nq.
\end{array}
$$
Noting that
$$
x_{q'+nq}=(T_{q'p'}\ybf)_{q'+nq}=y_{nq}+p'=y_0+np+p'=x_0+np+p',
$$
we then can make the configuration $\xbf$ by translating the segment $x_0,\ldots,x_{q'+nq}$ repeatedly by $T_{Q_n,P_n}$ and its inverse. In fact, $\xbf$ is just a $\delta$-gluing of  $\ybf$, $\zbf$, $T_{q'p'}\ybf$, and their translates by $T_{Q_nP_n}^m$, $m\in\Z$. Therefore, by continuity of $h_1$ and $h_2$, we have $|\dot{x}_j|<2\eps$ if we choose $\delta$ small enough.\emark

\begle\label{Fd}
Assume $\bar{F}>0$, $\xbf=(x_j)$ is a $(p,q)$-periodic configuration, and
$$
\dot{x}_j=-h_2(x_{j-1},x_j)-h_1(x_j,x_{j+1})+\bar{F}\leq 0\ \mbox{ for all }\ j\in\Z.
$$
Then $F_d(p/q)\geq \bar{F}$.
\ele
\noindent\pf
Assume on the contrary $0\leq F_d(p/q)<\bar{F}$. Then by Proposition~\ref{finite} there exists a $(p,q)$-periodic equilibrium $\ybf=(y_j)$ for\refer{tilt} with $F=F_d(p/q)$, and hence
$$
-h_2(y_{j-1},y_j)-h_1(y_j,y_{j+1})+\bar{F}> 0,
$$
implying together with the condition by Theorem 4.2 and Addendum 4.4 in~\cite{An} the existence of a $(p,q)$-periodic solution $\zbf=(z_j)$ of
$$
-h_2(z_{j-1},z_j)-h_1(z_j,z_{j+1})+\bar{F}= 0,\ j\in\Z,
$$
a contradiction to $F_d(p/q)<\bar{F}$ because of Proposition~\ref{finite}.\emark

\begde\label{family}
Let $\ybf$ and $\xbf$ be two Birkhoff $(p,q)$-periodic equilibria of\refer{tilt} and $\ybf\ll\xbf$. Assume there is a continuous mapping $\beta$ from an open interval $(a,b)$ to $\mathcal{X}$ such that

(i)\ $\beta(s_1)\ll\beta(s_2)$ for $s_1<s_2$, $s_1$, $s_2\in (a,b)$,

(ii)\ for each $s\in (a,b)$, $|\beta_n(s)-x_n|\to 0$ as $n\to +\infty$ and $|\beta_n(s)-y_n|\to 0$ as $n\to -\infty$, where $\beta_n$ is the projection of $\beta$ to site $n$,

(iii)\ for each $n\in\Z$, $\beta_n(s)\to x_n$ as $s\nearrow b$ and
$\beta_n(s)\to y_n$ as $s\searrow a$.

Then we say that $\mathcal{M}=\{\beta(s)\,|\,s\in (a,b)\}$ is a continuous family of advancing discommensurations between $\ybf$ and $\xbf$.
\ede

\begle\label{tangency}
Let $\ybf$ and $\xbf$ be two Birkhoff $(p,q)$-periodic equilibria of\refer{tilt} with $F=F_1$ and $\ybf\ll\xbf$. Assume there is a continuous family $\mathcal{M}$ of advancing discommensurations between $\ybf$ and $\xbf$, and for each $\zbf\in\mathcal{M}$, $\dot{\zbf}\geq \zero$ for\refer{tilt} with $F=F_1$. Suppose $\hat{\ybf}$ and $\hat{\xbf}$ are two Birkhoff $(p,q)$-periodic equilibria of\refer{tilt} with $F=F_2>F_1$ and
$\ybf\ll\hat{\ybf}<\xbf\ll\hat{\xbf}$. Then there is no equilibrium advancing discommensuration for\refer{tilt} with $F=F_2$ connecting $\hat{\ybf}$ to $\hat{\xbf}$.
\ele
\noindent\pf
Suppose there exists an equilibrium advancing discommensuration $\hat{\zbf}$ for\refer{tilt} with $F=F_2$ connecting $\hat{\ybf}$ to $\hat{\xbf}$, which implies
$$
|\hat{z}_n-\hat{x}_n|\to 0 \mbox{ as }n\to +\infty \mbox{ and }
|\hat{z}_n-\hat{y}_n|\to 0 \mbox{ as }n\to -\infty.
$$

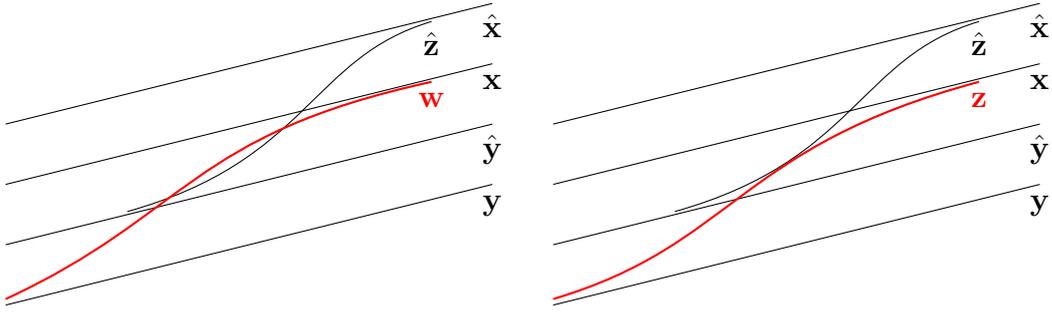
\begin{figure}
  \begin{center}






  \begin{tikzpicture}[scale=0.8]
  \draw (-9,0)--(-1,2)node [below]{$\ybf$};
  \draw (-9,1)--(-1,3)node [below]{$\hat{\ybf}$};
  \draw (-9,2)--(-1,4)node [below]{$\xbf$};
  \draw (-9,3)--(-1,5)node [below]{$\hat{\xbf}$};

  \draw (-7,1.55)..controls(-4,2.4) and (-4,4.1)..(-2,4.7)node [below]{$\hat{\zbf}$};
  \draw[red,thick] (-9,0.1)..controls(-6,1.5) and (-6,2.82)..(-2,3.7)node [below]{$\wbf$};


  \draw (0,0)--(8,2)node [below]{$\ybf$};
  \draw (0,1)--(8,3)node [below]{$\hat{\ybf}$};
  \draw (0,2)--(8,4)node [below]{$\xbf$};
  \draw (0,3)--(8,5)node [below]{$\hat{\xbf}$};

  \draw (2,1.55)..controls(5,2.4) and (5,4.1)..(7,4.7)node [below]{$\hat{\zbf}$};
  \draw[red,thick] (0,0.1)..controls(3,1) and (3,2.62)..(7,3.7)node [below]{$\zbf$};

  \end{tikzpicture}
  \end{center}
  \caption{The construction of $\zbf$.}\label{Fig:z}
  \end{figure}

Choose $s_0\in (a,b)$ such that $\wbf=(w_n)=\beta(s_0)\in\mathcal{M}$ and there exist $n_1<n_2<n_3$ with the following property (see Figure~\ref{Fig:z}),
$$
w_{n_1}<\hat{z}_{n_1},\ \  w_{n_2}\geq\hat{z}_{n_2}, \ \ \mbox{ and }\ w_{n_3}<\hat{z}_{n_3}.
$$
Let
$$
S=\{\beta(s)\in\mathcal{M}\,|\,\beta_n(s)\leq \hat{\zbf}_n,\ n_1\leq n\leq n_3,\ a<s\leq s_0\} \ \mbox{ and }\ \zbf=\sup S.
$$
Since $\beta_n(s)\to y_n$ uniformly for $n_1\leq n\leq n_3$ as $s\searrow a$,  the set
$S$ is nonempty and we have $\zbf\in\mathcal{M}$ due to the continuity of $\mathcal{M}$ (see Figure~\ref{Fig:z}). Moreover, from the property of $\wbf$, there exists $n\in\Z$, $n_1<n<n_3$,  such that $z_{n-1}\leq\hat{z}_{n-1}$, $z_n=\hat{z}_n$, and $z_{n+1}\leq \hat{z}_{n+1}$, i.e., $\zbf$ and $\hat{\zbf}$ are ``tangent'' at site $n$,  implying by the monotonicity
$$
\begin{array}{lll}
0\leq \dot{z}_n& = & -h_2(z_{n-1},z_n)-h_1(z_n,z_{n+1})+F_1\leq -h_2(\hat{z}_{n-1},\hat{z}_n)-h_1(\hat{z}_n,\hat{z}_{n+1})+F_1\\[4pt]
& < & -h_2(\hat{z}_{n-1},\hat{z}_n)-h_1(\hat{z}_n,\hat{z}_{n+1})+F_2,
\end{array}
$$
a contradiction to our assumption that  $\hat{\zbf}$ is an equilibrium of\refer{tilt} with $F=F_2$. Consequently, there is no equilibrium advancing discommensuration for\refer{tilt} with $F=F_2$ connecting $\hat{\ybf}$ to $\hat{\xbf}$.\emark

\subsection{Proof of Theorem A and its Corollary}
\label{sec:cty}
\noindent{\bf (a)\ Existence of limits.}\
Let $F_d(p/q+)={\lim\sup}_{\omega\searrow p/q}F_d(\omega)$. It suffices to show that
${\lim\inf}_{\omega\searrow p/q}F_d(\omega)=F_d(p/q+)$. Choose $\omega_n\searrow p/q$
as $n\to\infty$ such that $F_n=F_d(\omega_n)\to F_d(p/q+)>0$ as $n\to\infty$. Then from Theorem~\ref{construction-IOC-p/q} and its remark it follows that there exists an IOC $\ell$ of type $p/q+$ for\refer{tilt} with $F=F_d(p/q+)$ such that each gap of $M$ contains an equilibrium of type $p/q+$, where $M$ denotes the closed and strictly ordered set of equilibria on $\ell$ of type $(p,q)$.

Let $p'/q'$ be a Farey neighbour of $p/q$ with $p'/q' > p/q$ ({\em Farey neighbour} means $|p'q-pq'|=1$).

Let $\ybf=(y_n)$ be an equilibrium of type $(p,q)$ on $\ell$. Then $T_{q'p'}\ybf$ is also an equilibrium on $\ell$, which is above $\ybf$ and is the closest to $\ybf$ among the translations $\{T_{mn}\ybf\,|\,nq>mp\}$ above $\ybf$. Indeed, since $qp'-pq'=1$, we have
$$
\ybf\ll T_{q'p'}\ybf\ll T^2_{q'p'}\ybf\ll\cdots\ll T^q_{q'p'}\ybf=\ybf+\one.
$$

The simplest case is that
$\ybf$ and $T_{q'p'}\ybf$ are the endpoints of a gap in $M$, so let us assume this for the moment.
Then by hypothesis, there is an equilibrium $\zbf=(z_i)$ of type $p/q+$ connecting $\ybf$ to $T_{q'p'}\ybf$, i.e., $z_i-y_i\to 0$ as $i\to -\infty$, and $z_i-(T_{q'p'}\ybf)_i\to 0$ as $i\to +\infty$.
From Lemma~\ref{per-glue} we can construct, for each $\eps>0$,  a periodic configuration
$\xbf^n=(x_i^n)$ of type $(np+p',nq+q')$ such that $|\dot{x}_i^n|<2\eps$ for\refer{tilt} with
$F=F_d(p/q+)$, where $n\geq N$ for some $N$  depending on $\eps$.

If $\ybf$ and $T_{q'p'}\ybf$ are not endpoints of a gap in $M$, 
 then there is a non-empty countable (countable includes finite) set of gaps in between them, separated by equilibria of type $(p,q)$ or possibly by bands of equilibria of type $(p,q)$ (non-generic). We use equilibria of type $p/q+$ in each gap and a sequence of equilibria of type $(p,q)$ to obtain a configuration $\zbf=(z_i)$ in Lemma~\ref{per-glue} that is a $\delta$-gluing of these equilibria and satisfies $|\dot{z}_i|<\eps$ for $i\in\Z$ if $\delta$ is small enough. Then from Lemma~\ref{per-glue} we still have a configuration $\xbf^n=(x_i^n)$ of type $(np+p',nq+q')$ for $n\geq N$, such that $|\dot{x}_i^n|<2\eps$ for $i\in\Z$.

An example with one equilibrium between $\ybf$ and $T_{q'p'}\ybf$ is shown in
Figure~\ref{Fig:x}(b).

Let $m\in\N$ and $m\geq 2$. We remark that the construction of $\xbf^n$ for $n\geq N$ is the same as in Lemma~\ref{per-glue}. For fixed $n\geq N$, let
$$
\begin{array}{l}
x_i=x^n_i=y_i, \mbox{ for } i=0,1,\cdots,q'-1\ (\mbox{see the construction in Lemma~\ref{per-glue}}),\\[4pt]
x_i=x^n_i, \mbox{ for } i=q',\cdots,nq+q'-1,\\[4pt]
 x_i=x^N_i, \mbox{ for }i=nq+q',\cdots\cdots, nq+(m-1)Nq+mq'-1,\mbox{ and }\\[4pt]
x_i=(T_{q'p'}^m\ybf)_i, \mbox{ for }i=nq+(m-1)Nq+mq'.
\end{array}
$$
Then we make a configuration $\xbf=(x_i)$ with $|\dot{x}_i|< 2\eps$ for $i\in\Z$ as before by translating the segment $x_0,\cdots\cdots$,
$x_{nq+(m-1)Nq+mq'}$.

Now we have $|\dot{x}_i|< 2\eps$ and $\xbf$ is a periodic configuration of type $(P_n,Q_n)$, where
$P_n=np+(m-1)Np+mp'$ and $Q_n=nq+(m-1)Nq+mq'$, $n\geq N$. From Lemma~\ref{Fd} we have the conclusion that  for $m\in\N$ and $n\geq mN$,
$$
F_d(P/Q)\geq F_d(p/q+)-2\eps,\ \mbox{ in which }\ P=np+mp'\ \mbox{ and } \ Q=nq+mq'.
$$

Note that every rational $P/Q$ between $p/q$ and an upper Farey neighbour $p'/q'$ can be written as $P/Q = (np+mp')/(nq+mq')$ for some $n,m >0$. It is easy to check that for $m\in\N$ and $n<mN$,
$$
\frac{P}{Q}-\frac{p}{q}=\frac{np+mp'}{nq+mq'}-\frac{p}{q}>
\frac{Np+p'}{Nq+q'}-\frac{p}{q}=\frac{1}{q(Nq+q')}>0.
$$
Therefore, we arrive at the conclusion that for each $\eps>0$,
$$
{\lim\inf}_{P/Q\searrow p/q}F_d(P/Q)\geq F_d(p/q+)-2\eps,
$$
and hence
$$
{\lim\inf}_{\omega\searrow p/q}F_d(\omega)=F_d(p/q+)={\lim\sup}_{\omega\searrow p/q}F_d(\omega),
$$
since $F_d(\omega)$ is continuous at irrational points~\cite{QW1,WMWQ}.

From the proof of Theorem~\ref{construction-IOC-p/q} we know that for\refer{tilt} with tilt $F=F_d(p/q+)$ we have an IOC of type $p/q+$ on which there exists an equilibrium of type $(p,q)$, implying $F_d(p/q+)\leq F_d(p/q)$ by Proposition~\ref{finite}.



\vspace{2mm}


\noindent{\bf (b)\ Genericity.}\
To prove the genericity conclusions, we work in the $C^4$-topology on generating functions $h$, so $C^3$ for area-preserving twist maps.

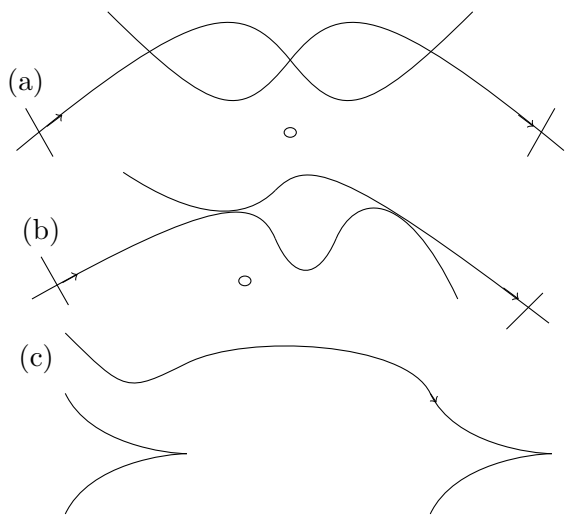
\begin{figure}[htbp] 
   \centering
\begin{tikzpicture} [scale=0.8]
\draw  (-4.5,-1.5)..controls(-1.5,1) and (-0.8,1)..(0,0);
\draw (0,0)..controls(0.8,-1) and (1.2,-1)..(3,0.8);
\draw (-3,0.8)..controls(-1.2,-1) and(-0.8,-1)..(0,0);
\draw (4.5,-1.5)..controls(1.5,1) and (0.8,1)..(0,0);
\draw (-3.9,-1.6)--(-4.35,-0.8) node[anchor=south]{\small (a)};
\draw (4.35,-0.8)--(3.9,-1.6);
\draw [->] (-4,-1.1)--(-3.75,-0.91);
\draw [->] (3.75,-0.91)--(4,-1.1);

\draw (0,-1.2)ellipse (0.1 and 0.08);
\end{tikzpicture}


\begin{tikzpicture}[scale=0.8]
\draw  (-4.5,-1.5)..controls(-1.5,0.2) and (-0.8,0.2)..(-0.5,-0.5);


\draw (-0.5,-0.5)..controls(-0.2,-1.2) and (0.2,-1.2)..(0.5,-0.5);

\draw (0.5,-0.5)..controls(0.8,0.2) and(1.6,0.43)..(2.5,-1.5);


\draw (-0.5,0.3)..controls(0,0.8) and (0.5,0.8)..(4,-1.9);
\draw (-0.5,0.3)..controls(-1,-0.2) and (-1.8,-0.2)..(-3,0.6);
\draw (-3.9,-1.6)--(-4.35,-0.8) node[anchor=south]{\small (b)};
\draw (3.3,-2)--(3.9,-1.4);
\draw [->] (-4,-1.23)--(-3.75,-1.09);
\draw [->] (3.25,-1.315)--(3.5,-1.5);

\draw (-1,-1.2)ellipse (0.1 and 0.08);
\end{tikzpicture}

\begin{tikzpicture}[scale=0.8]
\foreach \x in {0,6}
{
\draw  (-4+\x,0)..controls (-3.6+\x,-0.85)and (-2.3+\x,-1)..(-2+\x,-1);
\draw  (-4+\x,-2)..controls (-3.6+\x,-1.15) and (-2.3+\x,-1)..(-2+\x,-1);
}
\draw (-2,0.5)..controls(-3,0)..(-4,1) node[anchor=north east]{\small (c)};
\draw (-2,0.5)..controls(-1,1) and(1.6,0.85)..(2,0);

\draw [->] (2,0)--(2.1,-0.15);
\end{tikzpicture}

   \caption{Periodic orbits of type $(0,1)$ and their invariant manifolds for (a) $F=0$, (b) $F_d(0/1+)$ and (c) $F_d(0/1)$. 
   }
   \label{fig:manifolds}
\end{figure}

First we show that generically, $F_d(p/q+)>0$.  At $F=0$, by Aubry-Mather theory there is at least one minimising periodic orbit of type $(p,q)$ for the associated area-preserving map.  They are generically hyperbolic (use the formula for the trace of the derivative from \cite{MM}); in particular, there are finitely many modulo translations, but at least one and its translates.  For each gap in the set of minimising periodic orbits of type $(p,q)$ there is a minimising equilibrium advancing discommensuration from the left end to the right end (by Aubry-Mather theory again).  These are intersections of the stable and unstable manifolds if the periodic orbits are hyperbolic.  Generic intersections of the stable and unstable manifolds of hyperbolic periodic orbits are transverse (this follows from \cite{Rob}, because the twist condition is open).  Both the hyperbolic periodic orbits and the transverse intersections persist for small changes to $F>0$.  A typical picture is shown in Figure~\ref{fig:manifolds}(a) (for type $(0,1)$).

Let $p'/q'$ be the upper Farey neighbour of $p/q$.  First treat the simplest case when there is a single minimising periodic orbit of type $(p,q)$ (modulo translations).  Given $N>0$ large enough, take a configuration consisting of segments of length $N$ of homoclinic trajectory with small jumps near the periodic orbit.  They are ordered configurations of type $(Np+p',Nq+q')$, see the proof of Lemma~\ref{per-glue}.  Follow the gradient flow with tilt $F/2$ downhill to find a true periodic orbit.  Thus by Lemma~\ref{Fd}, for the given $F/2>0$ there is a sequence of equilibrium states with winding ratio converging from above to $p/q$, hence $F_d(p/q+)\ge F/2>0$. Actually, it is not necessary to construct ordered periodic orbits, because by the general theory of FK chains, if there is a periodic orbit of given type then there is an ordered one of that type.

If there is more than one gap, then one has to concatenate a sequence of heteroclinic trajectories crossing each gap in turn.

We next prove that generically $F_d(p/q+)<F_d(p/q)$.  At $F_d(p/q)$ there is generically a saddle-centre ordered periodic orbit of type $(p,q)$ and no other periodic orbits of type $(p,q)$, because there are none for $F>F_d(p/q)$.  It has stable and unstable semi-manifolds, as sketched in Figure~\ref{fig:manifolds}(c).  The semi-manifolds have no intersections and make at least one revolution.  There exists $\eps>0$ such that for $F \in (F_d(p/q)-\eps,F_d(p/q))$ the stable and unstable manifolds of the hyperbolic type $(p,q)$ orbit resulting from the saddle-centre continue to miss each other, at least in the first revolution.  If $F_d(p/q+) = F_d(p/q)$ then there exists for such $F$ a sequence of ordered equilibria of rotation numbers $\omega_n \searrow p/q$.  They have as a limit point an ordered equilibrium of rotation number $p/q+$ and it is asymptotic to ordered equilibria of type $(p,q)$.  Thus it has to be an ordered orbit of intersections of the stable and unstable manifolds of the hyperbolic periodic orbit of type $(p,q)$.  But this does not exist. 







\vspace{2mm}

\noindent {\bf (c)\ Proof of Corollary~\ref{co:conti}.}\
Continuity of $F_d$ at irrationals was already established in \cite{QW1}.  So it remains to prove continuity at $p/q\pm \in \mathcal{R}$.  One side is automatic from the definition, the other side is guaranteed by existence of the limit $\lim F_d(P/Q) = F_d(p/q+)$ as $P/Q \searrow p/q$, as shown before. This proves Corollary~\ref{co:conti}.


\vskip 1ex

{\noindent\bf Remark}: To interpret pictures in the phase space of the associated area-preserving twist map, like Figure~\ref{fig:manifolds},  generalisations of formulae in \cite{MMP} are useful.  See Appendix~C.


\subsection{Proof of Theorem B}\label{sec:pfB}
(i) Case $0<F\leq F_d(p/q+)$. From Theorem~\ref{construction-IOC-p/q} and its remark we know that there exists an IOC $\ell$ of type $p/q+$ for\refer{tilt} on which the set $M$ of equilibria of type $(p,q)$ is non-empty and closed, and in each gap of $M$ there is at least one equilibrium of type $p/q+$ whose translates by $T_{qp}^k$, $k\in\Z$ are also in the gap and accumulate pointwise on the endpoints of the gap.

Let $M'$ denote the union of $M$ with all the $p/q+$ equilibria on $\ell$. Then $M'$ is closed and $M'\ne \ell$ since $F>0$, as mentioned in the proof of Theorem~\ref{construction-IOC-p/q}. Note that an equilibrium of type $(p,q)$ is never an endpoint of a gap of $M'$ because it is accumulated on by type $p/q+$ equilibria. Consequently, each gap of $M'$ is invariant for\refer{tilt} with endpoints being equilibrium advancing discommensurations, and there exists a 
time-dependent solution $\xbf(t)$ in each gap such that $\xbf(t)$ tends to one end of the gap as $t\to -\infty$ and the other as $t\to +\infty$.  The direction depends on the action difference between the endpoints. 
\vskip 1ex
\noindent (ii) Case $F_d(p/q+)<F< F_d(p/q)$. Let $\eps=(F-F_d(p/q+))/2>0$. Choose $\omega_n\searrow p/q$ as $n\to\infty$ and $F>F_d(\omega_n)+\eps$. Let $\ell_n$ denote IOCs for\refer{tilt} with mean spacing $\omega_n$ and $\ell_F$ a limit point of $\ell_n$ in Hausdorff topology. Then $\ell_F$ is an IOC of\refer{tilt} with mean spacing $p/q$ and $T_{qp}\xbf \le \xbf$ for each $\xbf \in \ell_F$.

For each $\ybf=(y_i)\in\ell_n$, we have $\dot{y}_i>0$ for $i\in\Z$, by $F>F_d(\omega_n)$. So $\dot{x}_i\geq 0$, $i\in\Z$, for each $\xbf=(x_i)\in\ell_F$.

Also, $F\le F_d(p/q)$ implies that there is an equilibrium of type $(p,q)$ and hence an equilibrium $\ubf$ on $\ell_F$ of type $(p,q)$, and so each configuration on $\ell_F$ is bounded above and below by some translates of $\ubf$.

Let $M$ denote the set of equilibria of type $(p,q)$ on $\ell_F$. Then $M$ is nonempty, closed, and strictly ordered on $\ell_F$. Again $M\ne \ell_F$ since $F>0$.

It follows that $\ell_F$ consists of intervals $b$ of three types, each bounded by pairs $\xbf^-(b)\ll \xbf^+(b)$ of equilibria of type $(p,q)$:
\begin{enumerate}
\item type $p/q+$: configurations of type $p/q+$ left-asymptotic to $\xbf^-$, right-asymptotic to $\xbf^+$;
\item type $(p,q)$ sliding: configurations of type $(p,q)$ with $\dot{\xbf}\gg \zero$, backwards asymptotic to $\xbf^-$, forwards asymptotic to $\xbf^+$;
\item type $(p,q)$ equilibria: equilibria of type $(p,q)$.
\end{enumerate}
We call them {\em bands}.

From the remark of Theorem~\ref{construction-IOC-p/q} we know that there is at least one band of type $p/q+$.

If a band of type $p/q+$ contains no equilibria then $\dot{\zbf}\gg \zero$ for all its configurations. By the order property of the IOC and invariance under the gradient flow and $T_{qp}^{-1}$, there must be a $T>0$ at which $\phi^T_F(\zbf) =T_{qp}^{-1}\zbf$.  So we obtain a periodically sliding advancing discommensuration.

If there are no bands of type $(p,q)$ sliding, then 
 we deduce there is a band of type $p/q+$ containing no equilibria. Indeed,
suppose to the contrary that each gap of $M$ contains an equilibrium of type $p/q+$. Then similarly to the proof of part (a) of Theorem A we can construct a periodic configuration $\xbf$ of type $(np+p',nq+q')$, such that $|\dot{x}_i|<\eps$ for\refer{tilt}, $i\in\Z$, where $p'/q'>p/q$ is a Farey neighbour and $n\geq N$ for some $N\in\N$ depending on $\eps$. This implies from Lemma~\ref{Fd} that
$$
F_d(P_n/Q_n)\geq F-\eps= F_d(p/q+)+\eps,
$$
where $P_n=np+p'$, $Q_n=nq+q'$, and $P_n/Q_n\searrow p/q$ as $n\to\infty$, leading to a contradiction to Theorem A.
Thus there must be a band of type $p/q+$ containing no equilibrium of type $p/q+$ and hence a periodically sliding advancing discommensuration.

Next, we address what happens if all bands of type $p/q+$ contain an equilibrium, there is a band of type $(p,q)$ sliding and all type $(p,q)$ equilibria on $\ell_F$ are non-degenerate.  
We shall show that this could not happen owing to Theorem~\ref{positive-index}. Indeed, a band $b$ of type $(p,q)$ sliding corresponds to a trajectory of the gradient flow of\refer{tilt} in $\mcalX_{p,q}$ 
 backwards asymptotic to $\xbf^-(b)$ and forwards asymptotic to $\xbf^+(b)$, with $\dot{\xbf}\gg \zero$.  Thus $W_{p,q}$ decreases strictly along the trajectory and hence  $\xbf^-(b)$  has positive Morse index since it is non-degenerate. This is a contradiction to Theorem~\ref{positive-index}. As a consequence, we again have  a band of type $p/q+$ containing no equilibrium of type $p/q+$.

Lastly, we treat the case with a degenerate type $(p,q)$ equilibrium.  There is a band of type $p/q+$.  
 If it contains any equilibria then they are degenerate due to Theorem~\ref{degenerate-p/q+}.\emark

A periodically sliding discommensuration has an average velocity but in a different sense from that of solutions of type $(p,q)$ defined by\refer{nubar}.  Instead of considering $T$ for which $\xbf(T)=T_{01}\xbf(0)$ and defining $v=1/T$, for a periodically sliding advancing discommensuration we consider $T$ for which $\xbf(T) = T_{qp}^{-1}\xbf(0)$ and define $v=-1/T$.  Note that the sign convention makes $T>0$ (vacancies slide uphill). For a sliding retreating discommensuration we take $T$ for which $\xbf(T) = T_{qp}\xbf(0)$ and $v=1/T$.   These are velocities in $n$ rather than $x$.
\vskip1ex

An interesting possibility is the coexistence of periodically sliding and equilibrium advancing discommensurations. We make an example with $p/q = 0/1$ in Appendix B.






\subsection{Proof of Theorem C}
\label{sec:pfC}
(i) Assume $F_d(p/q+)=0$. Let $\omega_n\searrow p/q$ as $n\to\infty$. Let $\ell_n$ denote IOCs of mean spacing $\omega_n$ for\refer{tilt} with $F=F_d(\omega_n)\geq 0$. Then for each $\ybf=(y_i)\in\ell_n$, we have $\dot{y}_i\geq 0$ $\forall i\in \Z$, for\refer{tilt} with $F=F_d(\omega_n)$, see~\cite{QW1}. Let $\ell$ denote a limit point of $\ell_n$ in Hausdorff topology. Then for each $\xbf=(x_i)\in\ell$, $\dot{x}_i\geq 0$ for\refer{tilt} with $F=F_d(p/q+)=0$, i.e.,
\begin{equation}\label{untilt}
\dot{x}_i=-h_2(x_{i-1},x_i)-h_1(x_i,x_{i+1})\geq 0,\ \mbox{ for all }i\in\Z.
\end{equation}
Using the approach in the first part of the proof of Theorem A of~\cite{QW1}, we can furthermore prove that $\dot{x}_i=0$ for all $\xbf=(x_i)\in\ell$ and hence the existence of an invariant circle on the cylinder for the corresponding twist map. If each $\xbf\in\ell$ is of type $(p,q)$, then all orbits in the invariant circle are  $(p,q)$-periodic. If $\ell$ is type $p/q+$, then the invariant circle consists of some $(p,q)$-periodic orbits together with some orbits corresponding to equilibrium advancing discommensurations.
\vskip 1ex
(ii) In the other direction, assume for the twist map $\Phi_0$ ($F=0$ in\refer{eq:apmap}) there exists on the cylinder a rotational invariant circle consisting of some $(p,q)$-periodic orbits and some orbits corresponding to equilibrium advancing discommensurations.  By Birkhoff's theorem~\cite{Bi}, the circle is a graph.  Then we have an IOC $\ell$ for the corresponding untilted FK model\refer{tilt} with $F=0$, and each element on $\ell$ is an equilibrium. Let $\tau: \R\to {\R}^{\Z}$ be a continuous parametrisation of $\ell$, i.e., $\ell=\{\tau(s)\,|\,s\in\R\}$. In what follows we prove $F_d(p/q+)=0$ by contradiction.

Assume $F_d(p/q+)>0$. Then for each $0<\hat{F}\leq F_d(p/q+)$, from Theorem~\ref{construction-IOC-p/q} and its remark it follows that the tilted FK chain\refer{tilt}
with $F=\hat{F}$ has an IOC $\hat{\ell}$ of type $p/q+$ such that for any two neighbouring $(p,q)$-periodic equilibria $\hat{\ybf}\ll\hat{\xbf}$ on $\hat{\ell}$ there is an equilibrium advancing discommensuration $\hat{\zbf}$  of\refer{tilt} with $F=\hat{F}$ connecting $\hat{\ybf}$ to $\hat{\xbf}$.

We remark that any $(p,q)$-periodic equilibrium for\refer{tilt} with $F=0$ on $\ell$ and any $(p,q)$-periodic equilibrium for\refer{tilt} with $F=\hat{F}$ on $\hat{\ell}$ cannot cross (here ``cross'' means that their Aubry diagrams cross). Indeed, if $\xbf\in\ell$ and $\hat{\xbf}\in\hat{\ell}$ cross, say, $x_{n_1}<\hat{x}_{n_1}$, $x_{n_2}>\hat{x}_{n_2}$, $x_{n_3}<\hat{x}_{n_3}$, where $n_1<n_2<n_3$, then there exists in $\ell$ a configuration
$$
\wbf=\sup\{\ubf=(u_n)\in\ell\,|\,u_n\leq\hat{x}_n,\ n_1\leq n\leq n_3\},
$$
for which there is an $n\in\Z$, $n_1<n<n_3$, such that $w_{n-1}\leq \hat{x}_{n-1}$, $w_n=\hat{x}_n$, and $w_{n+1}\leq \hat{x}_{n+1}$ due to the continuity of $\tau$.
Note that because of the monotonicity we obtain
$$
0=-h_2(w_{n-1},w_n)-h_1(w_n,w_{n+1})\leq -h_2(\hat{x}_{n-1},\hat{x}_n)-h_1(\hat{x}_n,\hat{x}_{n+1}),
$$
implying
$$
-h_2(\hat{x}_{n-1},\hat{x}_n)-h_1(\hat{x}_n,\hat{x}_{n+1})+\hat{F}>0,
$$
 a contradiction to the fact that $\hat{\xbf}$ is an equilibrium of the tilted FK chain\refer{tilt} with $F=\hat{F}$.

Therefore, the union of equilibria of $(p,q)$ type on $\ell$ and those on $\hat{\ell}$ is an ordered set. Assume without loss of generality that between $\hat{\xbf}$ and $\hat{\ybf}$ there exists at least one $(p,q)$-periodic equilibrium on $\ell$ and let $\xbf\in\ell$ be the closest such one to $\hat{\ybf}$ (see Figure~\ref{Fig:z}), i.e.,
$$
\xbf=\inf\{\ubf\in\ell\,|,\hat{\ybf}\leq \ubf\leq \hat{\xbf}, \, \ubf\mbox{ is a }(p,q)\mbox{-periodic equilibrium}\}.
$$
Let
$$
\ybf=\sup\{\ubf\in\ell\,|\,\ubf\leq\hat{\ybf}, \,\ubf\mbox{ is a }(p,q)\mbox{-periodic equilibrium}\}.
$$
Then it follows that $\xbf\in\ell$, $\ybf\in\ell$, and $\ybf\ll\xbf$. Meanwhile, $\ybf$ and $\xbf$ are neighbouring equilibria of type $(p,q)$ on $\ell$, i.e., there are no $(p,q)$-periodic equilibria on $\ell$ in between $\ybf$ and $\xbf$.

From the assumption it follows that all configurations on $\ell$ between $\ybf$ and $\xbf$ are equilibrium advancing discommensurations for\refer{tilt} with $F=0$, that is, we have a continuous family $\mathcal{M}$ of advancing discommensurations between $\ybf$ and $\xbf$, and for all $\zbf\in\mathcal{M}$, $\dot{\zbf}=\zero$ for\refer{tilt} with $F=0$. Note that $\ybf\ll\hat{\ybf}<\xbf\ll\hat{\xbf}$. Then from Lemma~\ref{tangency} it follows that there is no equilibrium advancing discommensuration for\refer{tilt} with $F=\hat{F}$ connecting $\hat{\ybf}$ to $\hat{\xbf}$, a contradiction.
\emark




\newsection{Discussion}\label{sec:dis}

Our main result is existence of equilibrium discommensuration if $0< F < F_d(p/q\pm)$ and of periodically sliding discommensuration if $F_d(p/q\pm) < F < F_d(p/q)$.  So the simplest scenario is that there is a depinning transition for discommensurations as $F$ increases through $F_d(p/q\pm)$, analogous to that already established for spatially periodic configurations in \cite{BM}.

This result is likely to be important for understanding nonlinear conductivity in charge density wave materials~\cite{Gr,Th}.  If they have rational ratio $p/q$ of wavelengths of charge density to crystal structure, as is generically the case (though not necessarily of probability one), and the electric field $F$ is between $F_d(p/q\pm)$ and $F_d(p/q)$, then CB and RSM proposed the following scenario with Mario Floria back in 2005.  Thermal noise nucleates pairs of discommensurations, advancing and retreating.  They slide in opposite directions at speeds $v_\pm(F)$, producing a nett displacement of the charge density wave by one period (unless they get blocked by defects in the crystal lattice).  Thus conductivity is determined mostly by the nucleation rate of pairs of discommensurations.  This depends like $C\exp(-\beta\Delta W)$ on the energy barrier $\Delta W$ for creating the pair, which is a function of $F$, where $\beta$ is inverse temperature in energy units and $C$ is a relatively weak prefactor.
When $F>F_d(p/q)$ the charge density waves slide as a whole at speed $v(F)$ and the conductivity is proportional to this speed.  This effect is not likely to occur until significantly higher fields than $F_d(p/q\pm)$, however. So nucleation and movement of discommensurations in the regime $F_d(p/q\pm) < F < F_d(p/q)$ is likely to be the determining feature of the nonlinear conductivity of charge density wave materials.

Next, we make some technical comments.
\begin{enumerate}
\item Note that for $0 \le F < F_d(p/q)$, the ordered circle of type $(p,q)$ need not be unique.  An example is given in Appendix A.

\item The simplest scenario of depinning of discommensurations described above is not guaranteed:~it might be that some equilibrium discommensurations persist above $F_d(p/q\pm)$ and some sliding discommensurations exist already below $F_d(p/q\pm)$.  So there could be hysteresis.  Indeed,
in Appendix B we make  examples showing the coexistence of equilibrium discommensurations and periodically sliding discommensurations, for $F$ both above and below $F_d(0/1+)$. 

\item As remarked in Sec.~\ref{sec:gradflow} we have no need for the second derivative of $h$ to be bounded (see Appendix~E).

\item In Appendix~F, we give a useful representation for TOCs, in particular for IOCs.

\item We also note that the ordered circle approach gives analogues of the minimax sequences of Mather \cite{Mat2} (proposed also by \cite{ALA} but with an error in the construction, to be commented upon in Appendix A).  Namely, for $0\le F < F_d(p/q)$, in each gap between attracting equilibria in an ordered circle of type $(p,q)$ there is a repelling one.  The attracting ones generalise the minimising periodic sequences and the repelling ones generalise the minimax periodic sequences.  For $F=F_d(p/q)$ the equilibria are all attracting from the left and repelling to the right (else $F$ could be increased a little and equilibria would remain). The transition is topologically a saddle-node on invariant circle bifurcation.
\end{enumerate}


Finally, we pose some open questions.
\begin{enumerate}
\item For $F_d(p/q+)<F<F_d(p/q)$,  which states do the periodically sliding discommensurations attract?

\item  If there is an equilibrium state left- and right-asymptotic to different Birkhoff periodic ones, does there have to be a Birkhoff one (a property that holds for periodic states)?

\item We recall from \cite{BM} that when $F>F_d(p/q)$, all initial solutions of mean spacing $p/q$ are trapped between two phases of a unique periodically sliding solution of type $(p,q)$.  We suspect that there are IOCs of type $p/q\pm$ among these (it requires to show that the difference between the two phases does not go to zero).  Can this be proved?

\item Can one eliminate the exceptional case from Theorem B?
\end{enumerate}

\section*{Acknowledgements}
Qin is grateful to the support from the National Natural Science Foundation of China (11771316, 11790274, 12171347).
Baesens and MacKay are grateful to Soochow University for hospitality in Aug 2017.

\section*{Appendix A: Nonuniqueness of IOCs}
 \setcounter{theorem}{0}
 \setcounter{equation}{0}
 \setcounter{lemma}{0}
 \renewcommand{\theequation}{A.\arabic{equation}}
\renewcommand{\thelemma}{A.\arabic{lemma}}

Here we give an example of a generating function with more than one IOC of type $(1,2)$.  The example has zero tilt but one could tilt it a little and keep the non-uniqueness.  At the end we use the same example to discuss the construction of minimax orbits.

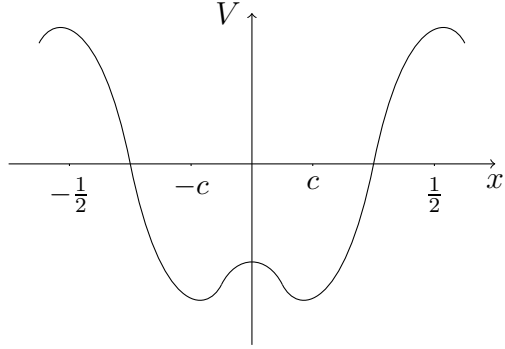
\begin{figure}[htbp] 
   \centering
\begin{tikzpicture}[scale=0.4]
  \draw[->] (-8,0)--(8,0)node [below]{$x$};
  \draw (-6,0)-- (-6,-2pt) node[below] {\small $-\frac{1}{2}$};
  \draw (-2,0)-- (-2,-2pt) node[below] {\small $-c$};
  \draw (2,0)-- (2,-2pt) node[below] {\small $c$};
  \draw (6,0)-- (6,-2pt) node[below] {\small $\frac{1}{2}$};
  \draw[->] (0,-6)--(0,5) node [anchor=east] {$V$};
  \draw (4,0)..controls(5,5) and (6.5,5)..(7,4);
  \draw (4,0)..controls(3,-5) and (1.5,-5)..(1,-4);
  \draw (1,-4)..controls(0.5,-3) and (-0.5,-3)..(-1,-4);
  \draw (-4,0)..controls(-3,-5) and (-1.5,-5)..(-1,-4);
  \draw (-4,0)..controls(-5,5) and (-6.5,5)..(-7,4);

  \end{tikzpicture}

   \caption{Potential $V$ for the example.}
   \label{fig:V}
\end{figure}

Let
$$h(x,x') = \frac{k}{2}(x'-x)^2 + V(x)$$ with
$$V(x)=-\frac{1}{4\pi^2}\cos{2\pi x} + \frac{b}{16\pi^2} \cos{4\pi x},$$
$b>1$ and $k>0$ small.
The even symmetry is not essential but makes the analysis simpler.
The potential $V$ has the following critical points, modulo $1$: maximum at $x=1/2$, local maximum at $x=0$, minima of the same height at $x = \pm c$ where $\cos{2\pi c} = 1/b$.  See Figure~\ref{fig:V}.

\begin{figure}[htbp] 
   \centering

\begin{tikzpicture}[scale=0.65]
\draw[->] (-1.5,0)--(8.5,0) node [anchor=west]{$x_0$};
\draw[->] (0,-0.5)--(0,9.5) node[anchor=east]{$x_1$};
\draw (-1.5,4)--(8.5,4);
\draw (8.5,8)--(-0.1,8) node [anchor=east] {$1$};
\draw (4,9.1)--(4,-0.1) node[below] {$1\over 2$};
\draw (8,9.1)--(8,-0.1)node[below] {$1$};
\draw[blue, thick,dashed] (-1.5,1)--(1,1)--(1,9)--(8.5,9);
\draw[red,  thick] (-1,0)--(-1,7)--(7,7)--(7,9.5);
\filldraw  (-1,1) circle (2pt);
\filldraw  (1,1) circle (2pt);
\filldraw  (1,7) circle (2pt);
\filldraw  (1,9) circle (2pt);
\filldraw  (7,9) circle (2pt);
\filldraw  (7,7) circle (2pt);
\filldraw  (-1,7) circle (2pt);

\draw [thick](-1.2,3.8)--(-0.8,4.2);
\draw [thick](-1.2,4.2)--(-0.8,3.8);
\draw [thick](-0.2,0.8)--(0.2,1.2);
\draw [thick](0.2,0.8)--(-0.2,1.2) node [anchor=east] {$c$};
\draw [thick](-0.2,6.8)--(0.2,7.2);
\draw [thick](0.2,6.8)--(-0.2,7.2);
\draw [thick](0.8,3.8)--(1.2,4.2);
\draw [thick](0.8,4.2)--(1.2,3.8);
\draw [thick](0.8,7.8)--(1.2,8.2);
\draw [thick](0.8,8.2)--(1.2,7.8);
\draw [thick](3.8,6.8)--(4.2,7.2);
\draw [thick](3.8,7.2)--(4.2,6.8);
\draw [thick](3.8,8.8)--(4.2,9.2);
\draw [thick](3.8,9.2)--(4.2,8.8);
\draw [thick](6.8,7.8)--(7.2,8.2);
\draw [thick](6.8,8.2)--(7.2,7.8);

\draw [thick] (-0.1,3.9) rectangle (0.1,4.1) node [anchor=south east] {$1\over 2$};
\draw [thick] (3.9,7.9) rectangle (4.1,8.1);

\filldraw (5,3) circle (2pt) node [anchor=west] {\small $\mbox{\ \ \ index\ }0$};
\draw [thick](4.8,1.8)--(5.2,2.2);
\draw [thick](4.8,2.2)--(5.0,2.0) node [anchor= west] {\small $\mbox{\ \ \ index\ }1$};
\draw [thick](5.0,2.0)--(5.2,1.8);
\draw [thick] (4.9,0.9) rectangle (5.1,1.1) node[anchor=west] {\small $\mbox{\ \ index\ }2$};

\draw (1,0)--(1,0.1) node [below=2pt] {$c$};
\draw (-1,0)--(-1,-0.1) node [below] {$-c$};
\draw (0,9)--(-0.1,9) node [anchor=east] {$1+c$};
\draw (0,7)--(-0.1,7) node [anchor=south east] {$1-c$};

\end{tikzpicture}

   \caption{Two IOCs in the plane of $(x_0,x_1)$, for the limiting case $k=0$.}
   \label{fig:2IOCs0}
\end{figure}
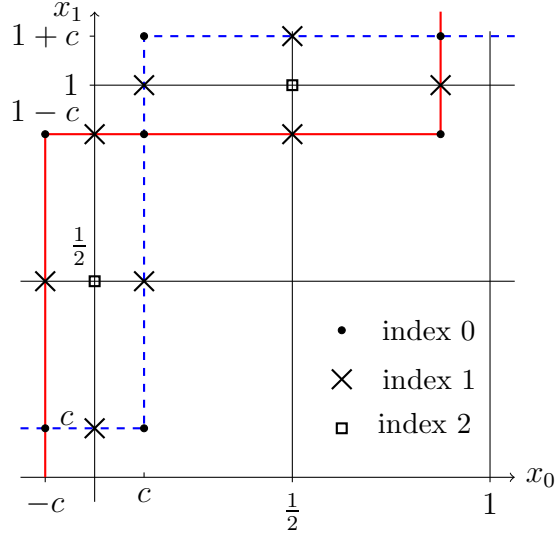

The energy for sequences of type $(1,2)$ is
\begin{equation}
W(x_0,x_1) = \tfrac{k}{2}(x_1-x_0)^2 + \tfrac{k}{2}(1+x_0-x_1)^2 + V(x_0) + V(x_1).
\end{equation}
Its gradient is
\begin{eqnarray}
\frac{\partial W}{\partial x_0} &=& k(1+ 2(x_0-x_1)) + V'(x_0), \\
\frac{\partial W}{\partial x_1} &=& k(-1+2(x_1-x_0)) + V'(x_1).
\end{eqnarray}
For the anti-integrable limit $k=0$ we form one IOC by connecting with horizontal and vertical lines the sequence of critical points
$$(-c,c), (0,c), (c,c), (c,\tfrac12), (c,1-c), (c,1), (c,(1+c), (\tfrac12, 1+c), (1-c,1+c).$$  Another IOC is given by connecting the sequence
$$(-c,c), (-c,\tfrac12), (-c,1-c), (0,1-c), (c, 1-c), (\tfrac12,1-c), (1-c,1-c), (1-c, 1), (1-c,1+c).$$
They are shown in Figure~\ref{fig:2IOCs0}.

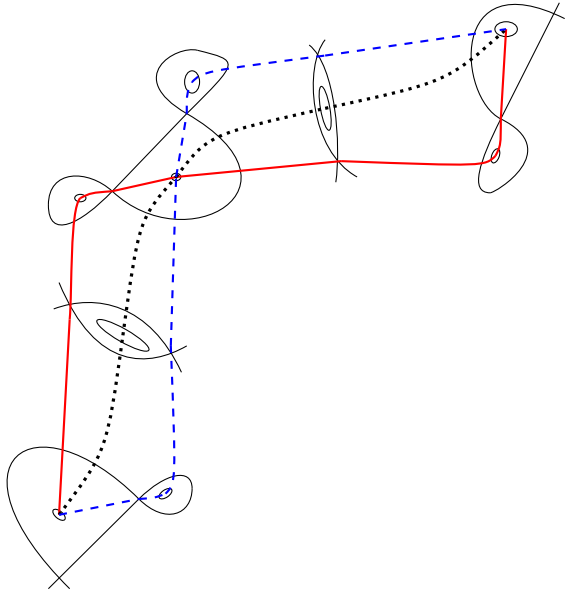
\begin{figure}[htbp] 
   \centering

\begin{tikzpicture}[scale=0.7]

\draw [rotate around={150:(0,4.2)}](0,4.2) ellipse (16pt and 3.5pt);
\draw (-1.2,5.2)..controls (-0.5,3.5) and (0.5,3.6)..(1.2,4);
\draw (-1.3,4.7)..controls (-0.5,5) and (0.5,4.8)..(1.1,3.5);

\draw (-0.8,6.8) ellipse (3pt and 2pt);
\draw (1,7.2) ellipse (2.5pt and 2pt);

\draw [rotate around={100:(3.8,8.5)}](3.8,8.5) ellipse (12pt and 2.5pt);
\draw (3.8,9.8)..controls (3.4,9.3) and (3.5,7.7)..(4.4,7.2);
\draw (3.5,9.7)..controls (4,9.3) and (4.1,7.8)..(4,7.1);

\draw [rotate around={90:(1.3,9)}](1.3,9) ellipse (6pt and 4pt);
\draw (-0.2,6.93)--(1.2,8.4);
\draw (-0.2,6.93)..controls(1.2,5.6) and (3.6,7)..(1.2,8.4);
\draw (-0.2,6.93)..controls(-0.66,7.35) and (-1.4,7.35)..(-1.4,6.67);
\draw (-0.2,6.93)..controls(-0.8,6.3) and (-1.4,6)..(-1.4,6.67);
\draw (1.2,8.4)..controls(0,9.1) and (1,9.7)..(1.4,9.6);
\draw (1.2,8.4)..controls(2.2,9.4)..(1.4,9.6);

\draw [rotate around={70:(7,7.6)}](7,7.6)ellipse (4pt and 2pt);
\draw(7.2,10) ellipse (6pt and 4pt);
\draw (7.1,8.3)..controls(6.3,8.7) and (6.1,11.3)..(8.3,10.2);
\draw (7.1,8.3)..controls(8.3,7.7) and (7,6.9)..(6.9,6.9);
\draw (7.1,8.3)..controls(6.7,7.5) and (6.5,6.9)..(6.9,6.9);
\draw (7.1,8.3)--(8.2,10.5);

\draw [rotate around={140:(-1.2,0.8)}](-1.2,0.8) ellipse (4pt and 2pt);
\draw [rotate around={35:(0.8,1.2)}](0.8,1.2) ellipse (4pt and 2pt);
\draw (0.3,1.1)..controls(-1.7,3.1) and (-3.4,1.8)..(-1,-0.6);
\draw (0.3,1.1)..controls(0.8,1.6) and (1.3,1.7)..(1.3,1.3);
\draw (0.3,1.1)..controls(0.8,0.6) and (1.3,0.8)..(1.3,1.3);
\draw (-1.4,-0.6)--(0.3,1.1);

\draw [red,  thick](-1.2,0.8)--(-1,4.5);
\draw [red,  thick](-1,4.5)..controls (-0.95,5.05)and (-1,6.6)..(-0.8,6.8);
\draw [red,  thick](-0.8,6.8)..controls (-0.6,6.9)..(-0.2,6.93);
\draw [red,  thick](-0.2,6.93)--(1,7.2);
\draw [red,  thick](1,7.2)--(4.1,7.5);
\draw [red,  thick](4.1,7.5)..controls (6.8,7.4)..(7,7.6);
\draw [red,  thick](7,7.6)..controls (7.1,7.8)..(7.1,8.3);
\draw [red,  thick](7.1,8.3)--(7.2,10);

\draw [dotted, very  thick] (-1.2,0.8)..controls (-0.3,2)..(0,4.2);
\draw [dotted, very thick] (0,4.2)..controls (0.3,6.4)..(1,7.2);
\draw [dotted, very thick] (1,7.2)..controls (1.6,8)..(3.8,8.5);
\draw [dotted, very thick] (3.8,8.5)..controls (6.1,9)..(7.2,10);

\draw [dashed,thick,blue](-1.2,0.8)--(0.3,1.1);
\draw [dashed,thick,blue] (0.3,1.1)..controls (0.6,1.12)..(0.8,1.2);
\draw [dashed,thick,blue] (0.8,1.2)..controls (1,1.4)..(0.9,4);
\draw [dashed,thick,blue](0.9,4)--(1,7.2);
\draw [dashed,thick,blue](1,7.2)--(1.2,8.4);
\draw [dashed,thick,blue] (1.2,8.4)..controls (1.2,8.8)..(1.3,9);
\draw [dashed,thick,blue] (1.3,9)..controls (1.5,9.2)..(3.8,9.5);
\draw [dashed,thick,blue](3.8,9.5)--(7.2,10);
\end{tikzpicture}

   \caption{Sketch of part of energy landscape for type $(1,2)$ for $k$ small positive, showing three IOCs.}
   \label{fig:2IOCs+}
\end{figure}

For $k>0$ small, the critical points persist and one obtains IOCs by following the gradient flow downhill from the chosen index 1 points.  This is illustrated in Figure~\ref{fig:2IOCs+} (confirmed by numerical computation for $K=0.03, b=2$ by Warwick project student Yifei Painter).

To prove that the resulting curves are IOCs, there are two key things to check:~firstly, that the critical points are ordered orbits, and secondly that the downhill eigenvector at the index $1$ points has positive slope.  After that, monotonicity of the gradient flow makes the resulting curves monotone.

The critical points are ordered orbits because they all start weakly ordered at $k=0$ and those which don't start strongly ordered use only local minima of $V$.  On making $k>0$ these become strongly ordered.

The eigenvectors at the index 1 points can be computed as follows.  The second derivative of $W$ is
\begin{equation}
D^2W(x_0,x_1) = \left[\begin{array}{cc}
2k+V''(x_0) & -2k \\
-2k & 2k + V''(x_1)
\end{array} \right] .
\end{equation}
So the eigenvalues are
$$\lambda = 2k + \frac12 (V''(x_0)+V''(x_1)) \pm \sqrt{\frac14 (V''(x_0)-V''(x_1))^2 + 4k^2}.$$
For index 1, the negative eigenvalue is obtained by taking the minus sign.  Its eigenvectors $(\xi_0,\xi_1)$ have slope
$$\xi_1/\xi_0 = \frac{1}{4k} (V''(x_0)-V''(x_1)) + \frac{1}{4k} \sqrt{(V''(x_0)-V''(x_1))^2+16k^2},$$
which is positive.

Actually, this example has a third IOC of type $(1,2)$.  Starting in the fast directions at the index 2 points $(0,1/2)$ and $(1/2, 1)$, one can follow the gradient flow to the minima of $W$ and it makes another IOC.  This is because the fast direction is given by taking the minus sign in the formula for the eigenvalue and the same formula for its eigenvector results so it has positive slope again.

Also, one can make examples with fewer periodic orbits.  Numerically, the critical points at $(c,c)$ and $(c,\tfrac12)$ annihilate each other as $k$ is increased, and the same for the other three similar pairs.  Taking gradient curves downhill from the remaining saddles produces two IOCs (and as before, there is a third one from the index 2 point).

To end this appendix, we use the same example to illustrate caveats in the construction of minimax orbits.

Given two minimisers of the same type $(p,q)$, it was proposed in \cite{ALA} to obtain a minimax state as follows: for each $x_0$ between those for the two minimisers take a minimum of the energy over $x_1,\ldots x_{q-1}$.  Then maximise the resulting energy over $x_0$.  Unfortunately, this procedure does not always produce a critical point of the energy, as was pointed out to RSM by Mather back in 1982.

\begin{figure}[htbp] 
   \centering
\begin{tikzpicture}[scale=0.7]
\draw (-1.9,3.1) parabola bend (0,6) (1.9,3.1);
\draw (-1.9,3.1) parabola bend (-0.6,5) (0,4.65);
\draw (0,4.65) parabola bend (0.6,5)(1.9,3.1);
\draw [very thick] (0,4.65)..controls (0.6,4)..(1.8,0);
\draw [very thick] (0,4.65)..controls (-0.6,4)..(-1.8,0);
\draw [->](-3,0)--(3,0) node [anchor=north]{$x_0$};
\draw [->](0,-0.4)--(0,7.2) node [anchor=east] {$E$};
\filldraw (-0.6,5) circle (2pt);
\filldraw (0.6,5) circle (2pt);
\filldraw (0,6) circle (2pt);
\draw [very thin] (1.3,5.5)--(1.301,5.5) node[anchor=south west]{\ \ \footnotesize index\ 2};
\draw [very thin](1.3,5.5)--(1.301,5.5) node[anchor=north west]{ \footnotesize minimax points};
\end{tikzpicture}

   \caption{Energy $E$ of critical points of $W$ with respect to $x_1$ as a function of $x_0$.  The thick line is the part obtained by minimising with respect to $x_1$.  The minimax points occur on its extension to local minima. }
   \label{fig:minimax}
\end{figure}
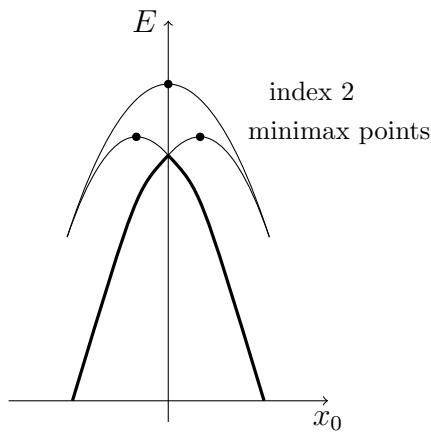

We illustrate this using the same example.  Looking at the level sets of $W_{1,2}$ in Figure~\ref{fig:2IOCs+} near $x_0=-0.5$ one sees that the line between the saddles has negative slope.
This can be proved by working out how the saddles $(\frac12,1\pm c)$ of the anti-integrable limit move for $k>0$.  The continuation of the saddles has $V'(x_0) = -k(1+2(x_0-x_1)) \sim \pm k c$.  For $x_0 = \frac12 + \xi$ then $V'(x_0) \sim -(1+b) \xi$.  Thus the saddles move to $x _0 \sim \frac12 \mp kc/(1+b)$.
Thus the minimum over $x_1$ is achieved at the tangency with a level set near the upper saddle for $x_0>-0.5$ and the lower saddle for $x_0 < -0.5$.  The maximum of this minimum occurs at $x_0=-0.5$ but does not correspond to a critical point.  If one includes all critical points of $W$ with respect to $x_1$ one obtains a diagram like Figure~\ref{fig:minimax}.

The desired minimax points are the two saddles.  They can be obtained by minimising the maximum energy over curves connecting the minima on either side.  An alternative is to take the infimum of $E$ such that both minima lie in the same connected component of the set $W(x_0,x_1)\le E$ and then argue that there must be a critical point on the boundary \cite{Mat2}.

The example has even symmetry which makes the saddles have the same height and so there are two minimax orbits.  One could break the even symmetry, however, to make one saddle lower than the other and hence a unique minimax orbit.

\section*{Appendix B: Coexistence of equilibrium and sliding discommensurations } 

 \setcounter{theorem}{0}
 \setcounter{equation}{0}
 \setcounter{lemma}{0}
 \renewcommand{\theequation}{B.\arabic{equation}}
\renewcommand{\thelemma}{B.\arabic{lemma}}


First we make an example with coexistence of equilibrium and sliding discommensurations of type $0/1+$ and show that it has $F$ above $F_d(0/1+)$.  Then we make a second example of coexistence and show that it has $F$ below $F_d(0/1+)$.

\begin{figure}[htbp] 
   \centering
\begin{tikzpicture}[scale=0.6]
\draw [->](-0.5,0)--(11,0)node [below]{$x$};
\draw [->](0,-0.5)--(0,5.5) node [anchor=west] {$V(x)-F x$};
\draw (1,0)--(-0.5,0) node [anchor=east] {$(a)$};
\draw (0,4)..controls (0.5,4)and (0.5,3)..(1,3);
\draw (1,3)..controls (1.5,3)and (1.5,5)..(2,5);
\draw (2,5)..controls (3,5)and (4,1)..(5,1);
\draw (5,1)..controls (5.5,1)and (5.5,2)..(6,2);
\draw (6,2)..controls (6.5,2)and (6.5,1)..(7,1);
\draw (7,1)..controls (7.5,1)and (7.5,3)..(8,3);
\draw (8,3)..controls (9,3)and (10,-1)..(11,-1);

\draw (1,3)--(1,0) node[below]{$a$};
\draw (5,1)--(5,0) node[below]{$b$};
\draw (7,1)--(7,0) node[below]{$a+1$};
\draw [dashed](2,5)--(2,0) node[below]{$c$};

\end{tikzpicture}

\begin{tikzpicture}[scale=0.6]
\draw (-0.5,0)--(1,0) node[below]{$a$};
\draw [->](1,0)--(11,0)node [below]{$X$};
\draw [->](0,-1.5)--(0,1.5) node [anchor=east]{$P$};
\draw (1,0)--(-0.5,0) node [anchor=east]{$(b)$};
\draw (3,0)..controls (3,0.6)and(1.5,1)..(0.8,-0.2);
\draw (0.8,0.9)..controls (1,0.4)and (2,1.2)..(3,1.2);
\draw (3,1.2)..controls (4,1.2)and (4.5,1)..(5,0) node [below]{$b$};
\draw (5,0)..controls (5.5,1)and (6.5,1)..(7.2,-0.2)node[below]{$a+1$};

\draw (0.8,0.2)..controls (1.5,-1)and (3,-0.6)..(3,0);
\draw (0.8,-0.9)..controls (1,-0.4)and (2,-1.2)..(3,-1.2);
\draw (3,-1.2)..controls (4,-1.2)and (4.5,-1)..(5,0);
\draw (5,0)..controls (5.5,-1)and (6.5,-1)..(7.2,0.2);

\draw (2.2,0) ellipse (0.2 and 0.1);
\draw (6,0) ellipse (0.2 and 0.1);

\draw[->] (5.9,0.71)--(6.1,0.71) node[above]{$X_1$};
\end{tikzpicture}


\begin{tikzpicture}[scale=0.6]
\draw (-0.5,0)--(1,0) node[below]{$a$};
\draw [->](1,0)--(11,0)node [below]{$X$};
\draw [->](0,-1.5)--(0,1.5) node [anchor=east] {$P$};
\draw (1,0)--(-0.5,0) node[anchor=east]{$(c)$};

\draw (0.8,0.2)..controls (1.5,-0.6)and(2.5,-0.4)..(2.5,0);
\draw (2.5,0) arc (0:180:0.4 and 0.2);
\draw (1.7,0) arc (180:360:0.3 and 0.15);
\draw (2.3,0) arc (0:180:0.2 and 0.1);

\draw (0.8,-0.2)..controls (2,1.5)and(4,1.5)..(5,0)node [below]{$b$};
\draw (5,0)..controls (5.5,-0.75)and(6.5,-0.5)..(6.5,0);
\draw (6.5,0) arc (0:180:0.25 and 0.13);
\draw (6,0) arc (180:360:0.125 and 0.063);
\draw (5.35,0) arc (180:360:0.5 and 0.25);
\draw (5.35,0)..controls (5.35,0.6)and(6.5,0.6)..(7.2,-0.2)node[below]{$a+1$};
\draw (5,0)..controls (5.5,1)and(6.8,0.2)..(7.2,0.6);
\draw (5,0)..controls (3,-1.8)and(1.5,-0.2)..(0.8,-0.6);

\draw (7.2,0.2)..controls (6.5,-1)and(5.5,-0.6)..(5,-0.6);
\draw (5,-0.6)..controls (4.5,-0.6)and(4,-1.2)..(3,-1.2);
\draw (3,-1.2)..controls (2,-1.2)and (1.5,-0.4)..(0.8,-1);

\draw[->] (3,1.1)--(3.2,1.1) node[above]{$X_2$};

\end{tikzpicture}


   \caption{(a) Example of $V(x)-Fx$;  (b) Connecting orbit $X_1$ for the steady states; (c) Connecting orbit $X_2$ for the 2D dynamics at speed $c$.}
   \label{fig:appB}
\end{figure}
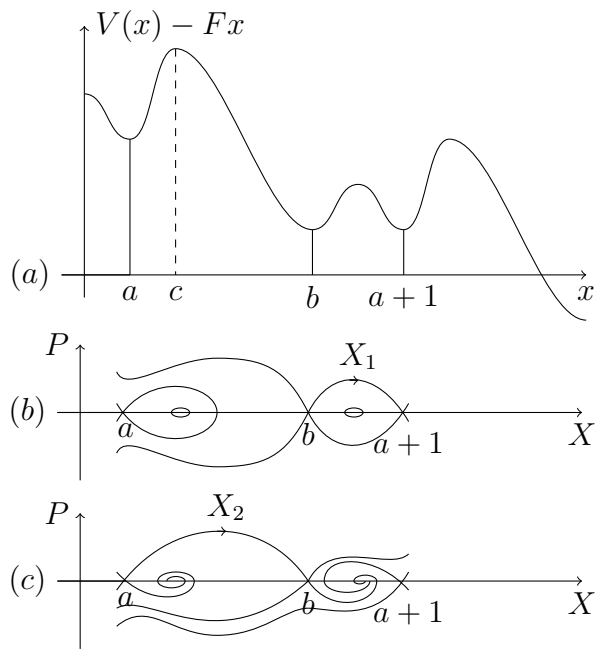

For the first example, we take
\begin{equation}
h(x,x') = \tfrac{k}{2} (x'-x)^2 + V(x)
\label{eq:hkV}
\end{equation}
with $k=\eps^{-2}$ large, i.e.~near the integrable (or continuum) limit.  We choose $V$ and $F>0$ so that $V(x)-Fx$ has two non-degenerate local minima per period.  Call their positions $a$ and $b$ with $a<b<a+1$.  We suppose $b$ and $a+1$ have the same value $E$ of $V(x)-Fx$.  See Figure~\ref{fig:appB}(a).

For the continuum limit with $s = n\eps$, the gradient flow 
of the energy
$$W = \int \frac12 \left(\frac{\partial x}{\partial s}\right)^2 + V(x)-Fx \ \mbox{d}s$$
is the 1D reaction-diffusion equation
\begin{equation}\label{diffusion}
\frac{\partial x}{\partial t} = \frac{\partial^2 x}{\partial s^2} - V'(x)+F.
\end{equation}
The reaction-diffusion equation has a steady solution 
 $x(s)=X_1(s+\tau)$ for any $\tau \in \R$, joining the unstable manifold of $b$ to the stable manifold of $a+1$, given by $P=\dot{X} = \sqrt{2(V(X)-FX-E)}$, as in Figure~\ref{fig:appB}(b).  $X_1$ is monotone.
By the standard theory of 1D reaction-diffusion equations, it also has a uniformly travelling front solution $x(t,s) = X_2(s+ct)$ for the value of $c>0$ such that the unstable manifold of $a$ joins the stable manifold of $b$, as in Figure~\ref{fig:appB}(c).  $X_2$ is monotone.

For small $\eps>0$, consider the discrete version of\refer{diffusion}
\begin{equation}\label{discrete}
\dot{x}_n=k(x_{n+1}-2x_n+x_{n-1})-V'(x_n)+F,
\end{equation}
where $x_n=x(t,n\eps)$, and $k=\eps^{-2}$.
Then the travelling front solution persists to one for\refer{discrete}, 
 near $x_n(t) = X_2(t+\eps c n)$, connecting the equilibria $a$ to $b$ \cite{CMPS, MP, BCC}.
That is, for large $k>0$, there exists a travelling front solution $\ybf(t)=(y_n(t))$, such that for each $t\in\R$, $y_n(t)<y_{n+1}(t)$ for $n\in\Z$, and $y_n(t)\to a$ as $t\to -\infty$,
$y_n(t)\to b$ as $t\to +\infty$, $y_n(t)\to a$ as $n\to -\infty$, $y_n(t)\to b$ as $n\to +\infty$, implying $\ybf(t)$ is a sliding advancing discommensuration for\refer{discrete}.

 Meanwhile, for each $k>0$, there exists~\cite{BCC} an equilibrium $\xbf=(x_n)$ such that
$x_n<x_{n+1}$ for $n\in\Z$, and $x_n\to b$ as $n\to -\infty$, $x_n\to a+1$ as $n\to +\infty$, implying $\xbf$ is an equilibrium advancing discommensuration of\refer{discrete}. 


{Therefore, for $k$ large enough, if we can show $F>F_d(0/1+, k)$, then we have the coexistence of equilibrium advancing discommensurations and periodically sliding advancing discommensurations above $F_d(0/1+)$, as desired.}



Assume to the contrary there exists a sequence $k_j\to +\infty$ as $j\to\infty$ such that $F\leq F_d(0/1+,k_j)$. We assume the inequality is strict for all $j\in\N$, the proof for the case that the equality holds for infinitely many $j$ being similar.

Fixing $k_j$, there exists an integer $q_j\geq 0$ such that $F<F_d(1/q,k_j)$ for $q>q_j$ and hence there is a Birkhoff periodic equilibrium $\xbf^q(k_j)$ of type $(1,q)$ for\refer{discrete} with $k=k_j$ . Then we construct $\xbf(k_j)$ as in the proof of case (i) of Theorem~\ref{construction-IOC-p/q}, replacing $\xbf^-$ and $\xbf^+$ there by $a\cdot\one$ and $c\cdot\one$ respectively, where $\one$ is the configuration with each component being $1$, and $c$ is the unique zero of $V'(x)-F$ between $a$ and $b$, see Figure~\ref{fig:appB}(a). Note that $\xbf(k_j)$ is actually the limit as $q\to\infty$ in product topology of a convergent subsequence of appropriate translates of the $\xbf^q(k_j)$.

Moreover, $\xbf(k_j)$ is a Birkhoff equilibrium with rotation number $0$ for\refer{discrete} with $k=k_j$ and $T_{10}\xbf(k_j)\leq \xbf(k_j)$. From the construction we know that there are two components of $\xbf(k_j)$, one $>a$  and the other $<c$, implying $\xbf(k_j)\neq a\cdot\one$ and $\xbf(k_j)\neq c\cdot\one$. As a consequence, $\xbf(k_j)$ is not a periodic equilibrium of type $(0,1)$ and hence $T_{10}\xbf(k_j)<\xbf(k_j)$, which means that $\xbf(k_j)$ is an equilibrium of type $0/1+$ for\refer{discrete} with $k=k_j$.
Let
$$
\xbf^-(k_j)=\lim_{m\to+\infty}T_{10}^m\xbf(k_j)\ \mbox{ and }\ \xbf^+(k_j)=\lim_{m\to-\infty}T_{10}^m\xbf(k_j).
$$
Then $\xbf^-(k_j)<\xbf^+(k_j)$ are two periodic equilibria of type $(0,1)$ for\refer{discrete}. Moreover, from  $0\leq x_n(k_j)-x_m(k_j)\leq 1$ for all $n>m$ we deduce that $\xbf^+(k_j)\leq \xbf^-(k_j)+\one$.

If we denote $\xbf^-(k_j)=a_1\cdot\one$ and $\xbf^+(k_j)=b_1\cdot\one$, then $a_1\leq a$ and $b_1\geq c$. Note that $b_1\neq c$ because the equilibrium $c\cdot\one$ has positive Morse index (see the proof of Theorem~\ref{positive-index}). Consequently, we have $b_1\geq b$. It then follows that $a_1=b-1$ or $a_1=a$, and $b_1=b$ or $b_1=a+1$. Without loss of generality we assume $a_1=a$, $b_1=b$, and there are infinitely many $k_j$, not relabelled, such that $\xbf^-(k_j)=a\cdot\one$ and $\xbf^+(k_j)=b\cdot\one$, implying $x_n(k_j)\to a$ as $n\to -\infty$ and $x_n(k_j)\to b$ as $n\to +\infty$, for each $j$.

There exists in product topology a convergent subsequence of $\{\xbf(k_j)\}$ since  $a\cdot\one\leq \xbf(k_j)\leq b\cdot\one$, not relabelled, such that 
$\zbf=(z_n)=\lim_{j\to+\infty}\xbf(k_j)$.
Note that $\xbf(k_j)=(x_n(k_j))$ is an equilibrium of\refer{discrete} with $k=k_j$, i.e.,
\begeq\label{kj}
k_j(x_{n+1}(k_j)-2x_n(k_j)+x_{n-1}(k_j))=V'(x_n(k_j))-F,\ \forall\,n\in\Z.
\eeq
Dividing by $k_j$ on both sides of\refer{kj} and letting $j\to\infty$, we have 
$$
z_{n+1}-z_n=z_n-z_{n-1}, \ \forall\,n\in\Z,  \ \mbox{ implying }\ z_{n+1}-z_n=0, \ \forall\, n\in\Z.
$$
As a consequence,
$$
0<\theta_n^j=x_{n+1}(k_j)-x_n(k_j)\to 0\ \mbox{ as }\ j\to\infty, \mbox{ for each }\ n\in\Z.
$$
We shall obtain a contradiction by multiplying  both sides of\refer{kj} by $\theta_n^j+\theta_{n-1}^j$.

Indeed, on one hand we have for each $k_j$,
$$
\sum_{n\in\Z}k_j(\theta_n^j-\theta_{n-1}^j)(\theta_n^j+\theta_{n-1}^j)=k_j\sum_{n\in\Z}(\theta_n^j)^2-(\theta_{n-1}^j)^2=0.
$$
On the other hand, it follows that
$$
\sum_{n\in\Z}\left(V'(x_n(k_j)-F\right)(\theta_n^j+\theta_{n-1}^j)\to 2\int_{a}^{b}V'(t)-F\mbox{d}t\neq 0,\ \mbox{ as }j\to\infty,
$$
a contradiction. Therefore, we conclude that $F>F_d(0/1+,k)$ if $k$ is large enough.

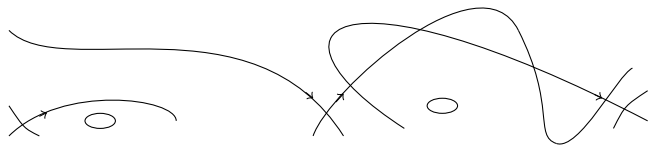
\begin{figure}[htbp] 
   \centering
\begin{tikzpicture}
\draw (0.8,1.6)..controls (1.5,2.3)and(3,2.1)..(3,1.8);
\draw (0.8,2)..controls (1,1.7)..(1.2,1.6);
\draw (2,1.8)ellipse (0.2 and 0.1);

\draw (0.8,3)..controls (1.5,2.3)and(4,3.5)..(5.2,1.6);
\draw (6,1.7)..controls (4,3)and(5,4)..(9.2,1.8);
\draw (4.8,1.6)..controls (5,2.2)and(7,4)..(7.5,3);
\draw (7.5,3)..controls (8,2)and(7.7,1.6)..(8,1.5);
\draw (8,1.5)..controls (8.3,1.4)and(8.8,2.4)..(9,2.5);
\draw (8.75,1.7)..controls (8.9,2)..(9.2,2.2);

\draw (6.5,2)ellipse (0.2 and 0.1);

\draw [->] (1.2,1.86)--(1.3,1.91);
\draw [->] (4.7,2.18)--(4.8,2.1);
\draw [->] (5.1,2.05)--(5.2,2.16);
\draw [->] (8.5,2.15)--(8.6,2.1);
\end{tikzpicture}
   \caption{Invariant manifolds for the two hyperbolic fixed points of the map.}
   \label{fig:invmflds}
\end{figure}

Note that generically, the invariant manifolds for the associated area-preserving map
\begin{eqnarray}
p' &=& p + \frac{1}{k}(V'(x)-F), \\
x' &=& x+p',
\end{eqnarray}
look like Figure~\ref{fig:invmflds}, so we can change $F$ a bit and keep the result.

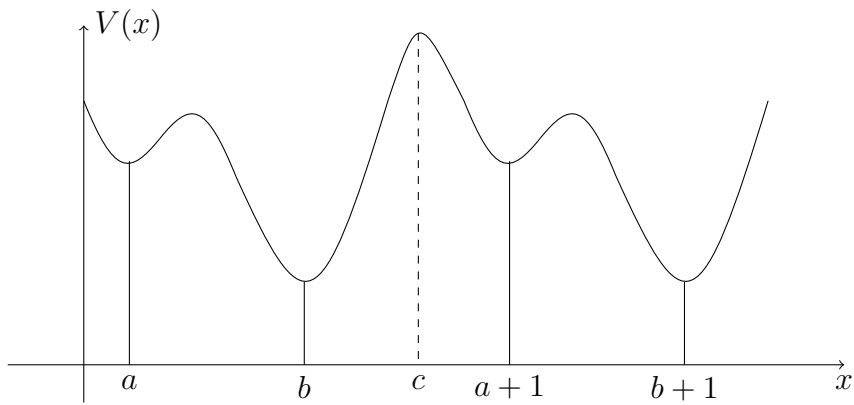
\begin{figure}[htbp] 
   \centering
\begin{tikzpicture}
\draw [->](-1,0)--(10,0)node [below]{$x$};
\draw [->](0,-0.5)--(0,4.5) node [anchor=west] {$V(x)$};

\draw (0,3.5)..controls (0.4,2.5)and(0.6,2.5)..(1,3);
\draw (1,3)..controls (1.4,3.5)and(1.6,3.5)..(2,2.5);
\draw (2,2.5)..controls (2.9,0.5)and(3.1,0.5)..(4,3.5);
\draw (4,3.5)..controls (4.4,4.7)..(5,3.5);

\draw (5,3.5)..controls (5.4,2.5)and(5.6,2.5)..(6,3);
\draw (6,3)..controls (6.4,3.5)and(6.6,3.5)..(7,2.5);
\draw (7,2.5)..controls (7.9,0.5)and(8.1,0.5)..(9,3.5);

\draw (0.6,2.7)--(0.6,0) node[below]{$a$};
\draw (2.9,1.1)--(2.9,0) node[below]{$b$};
\draw [dashed](4.4,4.4)--(4.4,0) node[below]{$c$};
\draw (5.6,2.7)--(5.6,0) node[below]{$a+1$};
\draw (7.9,1.1)--(7.9,0) node[below]{$b+1$};

\end{tikzpicture}

   \caption{Example of $V(x)$.}
   \label{fig:appB2}
\end{figure}

{Finally, we make an example showing the coexistence of equilibrium advancing discommensurations and periodically sliding advancing discommensurations below $F_d(0/1+)$.}

Take $h$ of the form (\ref{eq:hkV}) again, but choose
$V$  as in Figure~\ref{fig:appB2}. For $F=0$, the conclusions in~\cite{BCC} imply that for $k$ large, there is a sliding advancing discommensuration $\ybf(t)$ connecting $a$ and $b$ for\refer{discrete}. Moreover, there exists an equilibrium advancing discommensuration connecting $b$ and $b+1$, according to the Aubry-Mather theory. We need to verify that $F=0<F_d(0/1+)$. To this end, we require to construct the function $V$ to have $V''(c)<0$ large enough that the conclusion of \cite{BCC} applies for some $k$ satisfying $V''(c)<-2k$ (see Figure~\ref{fig:appB2}).   This can be done because although the size of $k$ in \cite{BCC} depends on the size of $V$ it does not depend on the size of $V''$.
Then following Criterion $1' $ in~\cite{MacP} (see also the approach used  in Appendix D), there are no rotational invariant circles crossing the vertical line $x=c$, {and hence $F_d(0/1+)>0$ by Theorem~C.}

\section*{Appendix C:~Interpretation of actions as areas}
 \setcounter{theorem}{0}
 \setcounter{equation}{0}
 \setcounter{lemma}{0}
 \renewcommand{\theequation}{C.\arabic{equation}}
\renewcommand{\thelemma}{C.\arabic{lemma}}

We generalise to the case of tilted FK chains some formulae of \cite{MMP} relating the areas under curves for an area-preserving twist map to the differences of action of pairs of orbits.

Recall that given a generating function $h$, which in this appendix includes the effect of tilt, we define an area-preserving map $\Phi: (x,y) \mapsto (x',y')$ (with notation change from $p$ to $y$) implicitly by
\begin{eqnarray}
y &=& -h_1(x,x') \label{eq:y} \\
y' &=& h_2(x,x'). \label{eq:y'}
\end{eqnarray}

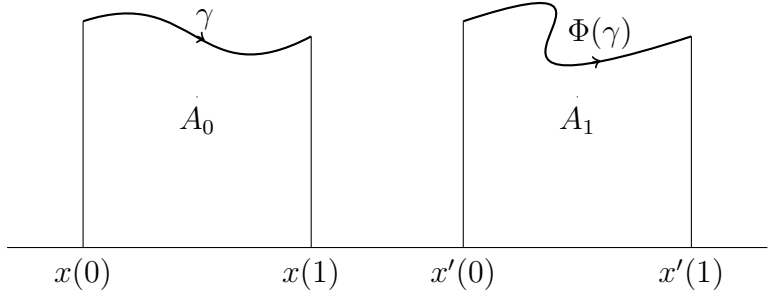
\begin{figure}[htbp] 
   \centering
   \begin{tikzpicture} 
\draw  (-1,0)--(9,0);
\draw  (0,3)--(0,0)node[below]{$x(0)$};
\draw  (3,2.8)--(3,0)node[below]{$x(1)$};
\draw  (5,3)--(5,0)node[below]{$x'(0)$};
\draw  (8,2.8)--(8,0)node[below]{$x'(1)$};
\filldraw (1.5,2)circle(0.05pt) node [below]{$A_0$};
\filldraw (6.5,2)circle(0.05pt) node [below]{$A_1$};
\draw [thick][->] (1.5,2.8)--(1.6,2.75) node [above]{$\gamma$};
\draw [thick][->] (6.7,2.45)--(6.8,2.48) node [above]{$\Phi(\gamma)$};

\draw [thick](0,3)..controls(1.5,3.5) and (1.5,2)..(3,2.8);
\draw [thick](5,3)..controls(8,4) and (4,1.5)..(8,2.8);

\end{tikzpicture}
   \caption{The areas $A_0$ and $A_1$ under $\gamma$ and $\Phi(\gamma)$ respectively.}
   \label{fig:A0A1}
\end{figure}

Given a differentiable arc $\gamma$ in $(x,y)$, parametrised by $s \in [0,1]$, the (oriented) area under it is defined to be $A_0 = \int_0^1 y \frac{\mbox{\scriptsize d}x}{\mbox{\scriptsize d}s} \mbox{d}s$.  Using (\ref{eq:y}), $A_0 = - \int_0^1 h_1(x,x') \frac{\mbox{\scriptsize d}x}{\mbox{\scriptsize d}s} \mbox{d}s$.  Then
 $$
\frac{\mbox{d}}{\mbox{d}s}h(x,x') = h_1(x,x')\frac{\mbox{d}x}{\mbox{d}s} + h_2(x,x')\frac{\mbox{d}x'}{\mbox{d}s}
$$
shows that
\begin{equation}
A_0 = h(x(0),x'(0))-h(x(1),x'(1)) + A_1
\label{eq:A01}
\end{equation}
where
$A_1 = \int_0^1 h_2(x,x') \frac{\mbox{\scriptsize d}x'}{\mbox{\scriptsize d}s} \mbox{d}s.$
But using (\ref{eq:y'}), the latter integral is the area under the image $\Phi(\gamma)$.
See Figure~\ref{fig:A0A1}.

If $\gamma$ is an arc of stable manifold, then the area $A_n$ under $\Phi^n(\gamma)$  which is represented by $(x_n(s),y_n(s))$, $s\in[0,1]$, goes to zero as $n \to \infty$, so iterating (\ref{eq:A01}) produces
\begin{equation}
A_0 = \sum_{n\ge 0} h(x_n(0),x_{n+1}(0))-h(x_n(1),x_{n+1}(1)).
\label{eq:Ws}
\end{equation}

Similarly, if $\tilde{\gamma}$ is an arc of unstable manifold from $(\tilde{x}(0),\tilde{y}(0))$ to $(\tilde{x}(1),\tilde{y}(1))$ then the area under it is $\sum_{n<0} h(\tilde{x}_n(1),\tilde{x}_{n+1}(1))-h(\tilde{x}_n(0),\tilde{x}_{n+1}(0))$.

In particular, if $(x(0),y(0))$ and $(x(1),y(1))$ are intersections of an unstable manifold and a stable manifold then the oriented area between the arcs of unstable and stable manifold connecting them is the difference in action of their orbits:
\begin{equation}
\Delta W = \sum_{n\in \Z} h(x_n(1),x_{n+1}(1))-h(x_n(0),x_{n+1}(0)).
\label{eq:DW}
\end{equation}

\begin{figure}[htbp] 
   \centering

   \begin{tikzpicture} 
\draw  [->](-0.5,0)--(10.5,0);
\draw  (0,3.25)--(0,0)node[below]{\footnotesize $x_0^-$};
\draw  (2,4.03)--(2,0)node[below]{\footnotesize$x_{-2q}+2p$};
\draw  (3.5,4.2)--(3.5,0)node[below]{\footnotesize$x_{-q}+p$};
\draw  (5,4.5)--(5,0)node[below]{\footnotesize$x_0$};
\draw  (6.5,4)--(6.5,0)node[below]{\footnotesize$x_q-p$};
\draw  (8,3.83)--(8,0)node[below]{\footnotesize$x_{2q}-2p$};
\draw  (10,3.05)--(10,0)node[below]{\footnotesize$x_0^+$};

\filldraw (1,1.5)circle(1pt) ++(0.2,0)circle(1pt) ++(0.2,0)circle(1pt);
\filldraw (2.75,2)circle(0.05pt) node [below]{$A_{-2}$};
\filldraw (4.25,2)circle(0.05pt) node [below]{$A_{-1}$};
\filldraw (5.75,2)circle(0.05pt) node [below]{$A_0$};
\filldraw (7.25,2)circle(0.05pt) node [below]{$A_1$};
\filldraw (9,1.5)circle(1pt) ++(0.2,0)circle(1pt) ++(0.2,0)circle(1pt);

\draw [thick][->] (1.2,3.78)--(1.3,3.82) node [above]{$W^u$};
\draw [thick][->] (8.8,3.58)--(8.9,3.52) node [above]{$W^s$};

\draw [thick](-0.5,3)..controls(2,4.2)..(5,4.2);
\draw [thick](5,4)..controls(8,4)..(10.5,2.8);
\draw [thick](-0.3,3.7)--(0.3,2.8);
\draw [thick](9.7,2.6)--(10.3,3.5);

\end{tikzpicture}

   \caption{The areas under pieces of stable and unstable manifold of periodic points of type $(p,q)$.}
   \label{fig:areaWsu}
\end{figure}
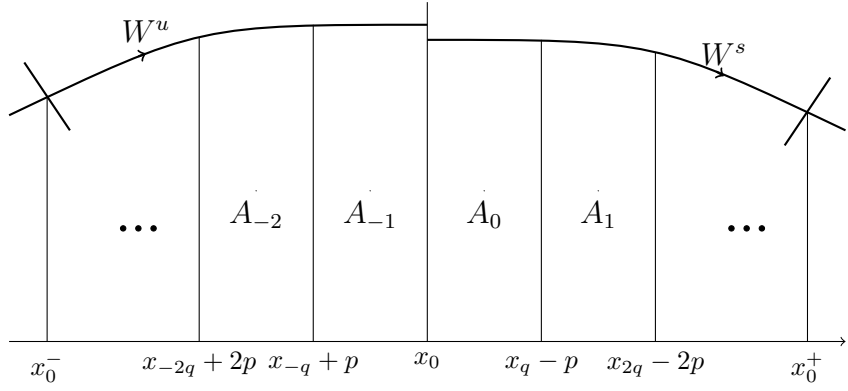

In the case that $\gamma$ is an arc of stable manifold of a type $(p,q)$ orbit
  $\xbf^+=(x_n^+)$  from a point $(x_0,y_0)$ to $\Phi^q(x_0,y_0)-(p,0)$, we can instead iterate (\ref{eq:A01}) $q$ times to obtain
$$A_0 = \sum_{n=0}^{q-1} [h(x_n,x_{n+1})-h(x_{q+n}-p,x_{q+n+1}-p)] + A_q$$
and then iterate this $k$ times to obtain a telescoping sum:
$$A_0 = \sum_{n=0}^{q-1} [h(x_n,x_{n+1})-h(x_{kq+n}-kp,x_{kq+n+1}-kp)]+A_{kq}.$$
As $k \to +\infty$, $A_{kq}\to 0$ and $x_{kq+n}-kp \to x^+_n$.  Thus
$$A_0 = \sum_{n=0}^{q-1} h(x_n,x_{n+1})-h(x^+_{n},x^+_{n+1}).$$
Then the area under the stable manifold from $(x_0,y_0)$ to $(x^+_0,y^+_0)$ is
$$\sum_{k\ge 0} A_{kq} = \sum_{n=0}^\infty h(x_n,x_{n+1})-h(x^+_{n},x^+_{n+1}).$$

Similarly, the area under the unstable manifold of a periodic point $(x^-_0,y^-_0)$ to a point $(x_0,\tilde{y}_0)$ on the unstable manifold is
$$\sum_{n<0} h(x_n,x_{n+1})-h(x^-_{n},x^-_{n+1}).$$
See Figure~\ref{fig:areaWsu}.

Summing these two, we see that the area under a curve formed of unstable manifold from a periodic point $(x^-_0,y^-_0)$ to a point $(x_0,\tilde{y}_0)$  followed by a vertical jump to $(x_0,{y}_0)$ on the stable manifold of a periodic point $(x^+_0,y^+_0)$ and then along the stable manifold to $(x^+_0,y^+_0)$ is
$$W_{\xbf^-,\xbf^+}(\xbf) = \sum_{n<0} [h(x_n,x_{n+1})-h(x^-_{n},x^-_{n+1})] + \sum_{n\ge 0} [h(x_n,x_{n+1})-h(x^+_{n},x^+_{n+1})],$$
where $\xbf = (x_n)$ is the sequence formed by taking the backwards orbit of $(x_0,\tilde{y}_0)$ and the forwards orbit of $(x_0,{y}_0)$.

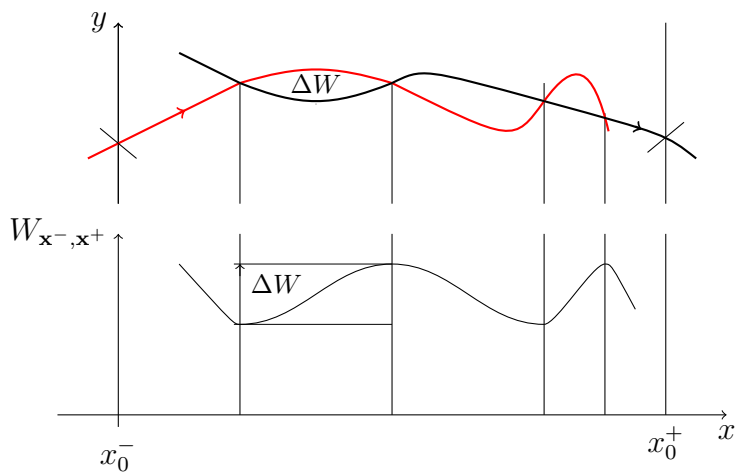
\begin{figure}[htbp] 
   \centering

  \begin{tikzpicture} [scale=0.8]
\draw  [->](-1,0)--(10,0)node[below]{$x$};
\draw [->](0,-0.2)node [below]{$x_0^-$}--(0,3)node [anchor=east]{$W_{\xbf^-,\xbf^+}$};
\draw [->] (0,3.5)--(0,6.5);
\draw [thin](2,0)--(2,3);
\draw [thin](4.5,0)--(4.5,3);
\draw [thin](7,0)--(7,3);
\draw [thin](8,0)--(8,3);
\draw [thin](9,0)node [below]{$x_0^+$}--(9,3);
\draw [very thin](1.9,1.5)--(4.5,1.5);
\draw [very thin](1.9,2.5)--(4.5,2.5);
\draw [->] (2,2.3)--(2,2.5)node [anchor=north west]{\footnotesize $\Delta W$};

\draw (1,2.5)..controls(1.9,1.5)..(2,1.5);
\draw (2,1.5)..controls(3,1.5) and (3.5,2.5)..(4.5,2.5);
\draw (4.5,2.5)..controls(5.5,2.5) and (6,1.5)..(7,1.5);
\draw (7,1.5)..controls(7.2,1.5) and (7.8,2.5)..(8,2.5);
\draw (8,2.5)..controls(8.1,2.5)..(8.5,1.75);

\draw [->](0,3.5)--(0,6.5)node[anchor=east]{$y$};
\draw [thin] (2,5.5)--(2,3.5);
\draw [thin] (4.5,5.5)--(4.5,3.5);
\draw [thin] (7,5.5)--(7,3.5);
\draw [thin] (8,5)--(8,3.5);
\draw [thin] (9,6.5)--(9,3.5);

\draw [thick, red](-0.5,4.25)..controls(1,5)..(2,5.5);
\draw [thick, red](2,5.5)..controls(3,5.8)and (3.5,5.8)..(4.5,5.5);
\draw [thick, red](4.5,5.5)..controls(6.5,4.5)..(7,5.2);
\draw [thick, red](7,5.2)..controls(7.5,5.9)and(7.8,5.8)..(8.06,4.7);

\draw [thick](2,5.5)..controls(3,5.1)and(3.5,5.1)..(4.5,5.5);
\draw [thick](1,6)--(2,5.5);
\draw [thick](4.5,5.5)..controls(5,5.75)..(7,5.2);
\draw [thick](7,5.2)..controls(9,4.65)..(9.5,4.25);

\draw (-0.3,4.75)--(0.3,4.25);
\draw (8.7,4.35)--(9.3,4.85);
\draw [thick,red][->](1,5)--(1.1,5.05);
\draw [thick][->](8.5,4.78)--(8.6,4.73);

\draw (3.25,5.15)circle(0.03pt)node[above]{\footnotesize$\Delta W$};

\end{tikzpicture}

   \caption{Illustration of variation of $W_{\xbf^-,\xbf^+}(\xbf)$ along a path formed by taking stable manifold forwards and unstable manifold backwards from $\xbf$.}
   \label{fig:DeltaW}
\end{figure}

In particular, if $(x_0,y_0)$ and $(\tilde{x}_0,\tilde{y}_0)$ are intersections of the unstable manifold of $\xbf^-$ and stable manifold of $\xbf^+$
then the area under the unstable manifold joining them minus the area under the stable manifold joining them equals the action difference
$W_{\xbf^-,\xbf^+}(\tilde{\xbf})-W_{\xbf^-,\xbf^+}(\xbf)$, which agrees with (\ref{eq:DW}).  The orbits of such intersections are critical points of $W_{\xbf^-,\xbf^+}$.  The height difference between the critical points is $\Delta W$.  This is illustrated in Figure~\ref{fig:DeltaW}.

Furthermore, if $\xbf=(x_n)$ is a local minimum of $W_{\xbf^-,\xbf^+}$ then at $(x_0,y_0)$ the tangent to the unstable manifold pointing away from $(x^-_0,y^-_0)$ is anticlockwise from the tangent to the stable manifold pointing towards $(x^+_0,y^+_0)$ because varying $\xbf$ can not decrease the area under the curve with a jump at $\xbf$.

So far in this appendix we have considered the area represented by the difference in  action between two orbits which converge together as  $n \to \infty$, and the area represented by the action of a pseudo-orbit that consists of a half-orbit converging to one periodic orbit of type $(p,q)$ as $n\to -\infty$ and to another as $n\to +\infty$ with a jump at $x_0$, and the special case where there is no jump.  To complete the appendix we provide a graphical interpretation of IOCs in the phase space for the area-preserving map.

Given an IOC we can associate two graphs $Y^-$ and $Y^+$ over $\R/\Z$ by
\begin{eqnarray}
Y^+(x_0) = -h_1(x_0,x_1), \nonumber \\
Y^-(x_0) = h_2(x_{-1},x_0), \nonumber
\end{eqnarray}
as $\xbf=(x_n)$ varies along the IOC.  By the order property of the IOC, $x_0$ varies continuously and strictly monotonically over the configurations of the IOC.  Also by the translation invariance, the same functions $Y^\pm$ are obtained whichever point is selected to be $n=0$.

\begin{figure}[htbp] 
   \centering

   \begin{tikzpicture}
\draw  [->](-0.5,0)--(8.5,0)node[below]{$x$};
\draw [->] (0,-0.5)--(0,5)node [anchor=east]{$y$};
\draw (8,0)--(8,5);
\draw [thick,red](0,2)..controls(3,2) and (3,3.5)..(4,3.5);
\draw [thick,red](4,3.5)..controls(5,3.5) and (5,2)..(8,2);
\draw [thick](0,4)..controls(3,4) and (3,2.3)..(4,2.3);
\draw [thick](4,2.3)..controls(5,2.3) and (5,4)..(8,4);

\draw (0.6,0.05)--(0.6,-0.05)node[below]{$x_{-1}$};
\draw (1.2,0.05)--(1.2,-0.05)node[below]{$x_0$};
\draw (1.8,0.05)--(1.8,-0.05)node[below]{$x_1$};

\filldraw (0.6,3.96)circle(1.2pt);
\filldraw (1.2,2.1)circle(1.2pt)node[below]{$Y^-$};
\filldraw (1.2,3.9)circle(1.2pt)node[above]{$Y^+$};
\filldraw (1.8,2.24)circle(1.2pt);

\filldraw (1.2,3.9)circle(1.2pt)node[anchor=north east]{$x_0$};

\draw[thin](0.6,3.96)--(1.2,2.1);
\draw[thin](1.2,2.1)--(1.2,3.9);
\draw[thin](1.2,3.9)--(1.8,2.24);
\draw [->](0.6,3.96)--(0.9,3.02)node [anchor=north east]{$\Phi$};
\draw[->](1.2,2.1)--(1.2,3.7);
\draw[->](1.2,3.9)--(1.5,3.07)node[anchor=south west]{$\Phi$};

\end{tikzpicture}

   \caption{An IOC corresponds to a pair of graphs of functions $(Y^-,Y^+)$ in the phase space for the area-preserving map $\Phi$.  The configurations on the IOC can be recovered by iteration.}
   \label{fig:IOCinxy}
\end{figure}
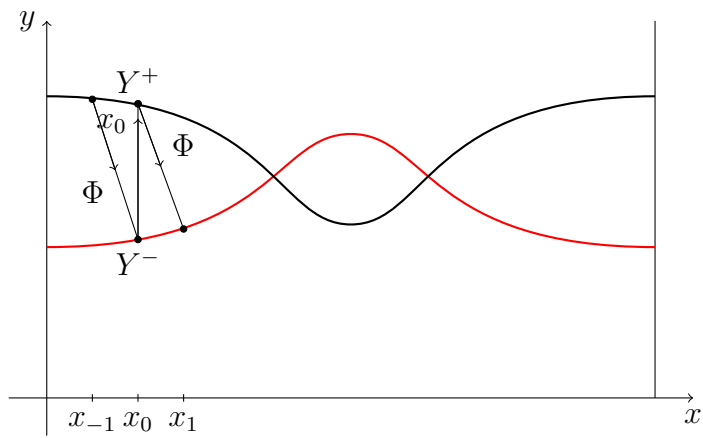

The difference between the graphs is precisely $\dot{x}_0$:
$$\dot{x}_0 = -h_2(x_{-1},x_0)-h_1(x_0,x_1)=Y^+(x_0)-Y^-(x_0).$$

The configurations of an IOC can be recovered from the pairs of graphs by iteration.  Given $x_0$, obtain $x_1$ by applying $\Phi$ to $(x_0,Y^+(x_0))$ and obtain $x_{-1}$ by applying $\Phi^{-1}$ to $(x_0,Y^-(x_0))$.  See Figure~\ref{fig:IOCinxy}.

\section*{Appendix D:~Example with a rotational invariant circle not of type $(p,q)$ nor $p/q+$ nor $p/q-$}
 \setcounter{theorem}{0}
 \setcounter{equation}{0}
 \setcounter{lemma}{0}
 \renewcommand{\theequation}{D.\arabic{equation}}
\renewcommand{\thelemma}{D.\arabic{lemma}}

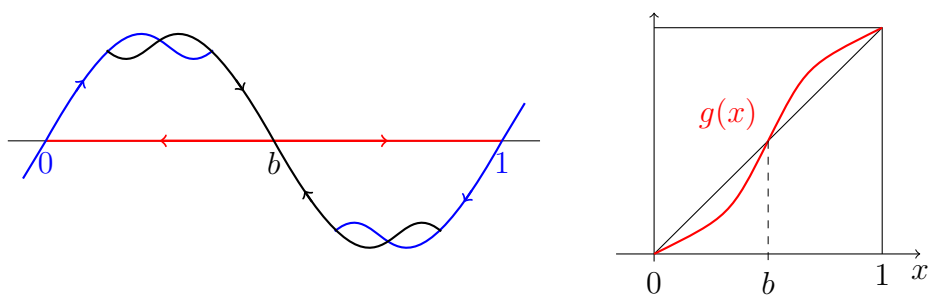
\begin{figure}[htbp] 
   \centering

   \begin{tikzpicture}
\draw (-0.5,0)--(6.5,0);
\draw[thick,red](0,0)--(6,0);
\draw [thick,red][->](3,0)--(1.5,0);
\draw [thick,red,->](3,0)--(4.5,0);
\draw [thick,blue](0,0)..controls(1.5,2.7)and(1.5,0.6)..(2.2,1.2);
\draw [thick,->,blue](0.4,0.66)--(0.5,0.81);

\draw [thick](3,0)..controls(1.5,2.7)and(1.5,0.6)..(0.8,1.2);
\draw [thick,->,](2.5,0.81)--(2.6,0.66);

\draw [thick,blue](6,0)..controls(4.5,-2.7)and(4.5,-0.6)..(3.8,-1.2);
\draw [thick,->,blue](5.6,-0.66)--(5.5,-0.81);

\draw [thick](3,0)..controls(4.5,-2.7)and(4.5,-0.6)..(5.2,-1.2);
\draw [thick,->,](3.5,-0.81)--(3.4,-0.66);

\draw [thick,blue](-0.3,-0.5)--(0,0)node [below]{$0$};
\draw [thick,blue](6.3,0.5)--(6,0)node[below]{$1$};

\filldraw(3,0)circle(0.05pt)node[below]{$b$};

\draw [->](7.5,-1.5)--(11.5,-1.5)node[below]{$x$};
\draw[->](8,-1.6)node[below]{$0$}--(8,1.7);
\draw (8,-1.5)--(11,1.5);
\draw (11,-1.5)node[below]{$1$}--(11,1.5)--(8,1.5);

\draw [red,thick](8,-1.5)..controls(9,-1)..(9.5,0);
\draw [red,thick](9.5,0)..controls(10,1)..(11,1.5);

\filldraw[red](9.5,0)circle(0.05pt)node[anchor=south east]{$g(x)$};
\draw [dashed,very thin](9.5,0)--(9.5,-1.6)node [below]{$b$};
\end{tikzpicture}

   \caption{(a) The hyperbolic fixed points of $\Phi$ and their invariant manifolds; (b) the circle diffeomorphism $g$. }
   \label{fig:funnyric}
\end{figure}

To illustrate the gap between Theorem C and a complete answer to the question of existence of rotational invariant circles for area-preserving twist maps with zero nett flux, we construct an example in which the only rotational invariant circle of rotation number $0$ consists of two hyperbolic fixed points together with a left-going and a right-going separatrix from one to the other.  See the horizontal red curve in Figure~\ref{fig:funnyric}(a).

We construct the example by a discrete analogue of Ma\~n\'e's Lagrangians.  Namely, we take $$h(x,x')= \frac12(x'-g(x))^2$$ for a lift $g$ of a degree-one circle diffeomorphism, chosen to have two fixed points $0$ and $0<b<1$ with $g'(0)\ne 1$, $g'(b)\ne 1$, and $g(x)<x$ for $x \in (0,b)$, $g(x)>x$ for $x \in (b,1)$, as in Figure~\ref{fig:funnyric}(b).  $h$ satisfies the conditions, notably $h_{12}(x,x') = -g'(x)<0$.

Every orbit of $g$ minimises the action sum, hence gives an orbit of the associated area-preserving twist map $\Phi$.  Or more directly, $\Phi$ is given implicitly by
\begin{eqnarray}
y &=& g'(x)(x'-g(x)), \nonumber \\
y' &=& x'-g(x), \nonumber
\end{eqnarray}
so these orbits all lie on $y=0$.
The fixed points of $g$ give hyperbolic fixed points of $\Phi$.  The hyperbolicity can be checked by finding the eigenvalues $\lambda,1/\lambda$ of $D\Phi$ there, or slightly faster by using the formula \cite{MM}
$$\lambda+1/\lambda-2 = \frac{\det D^2W_{p,q}}{\prod_{n=0}^{q-1} (-h_{12}(x_n,x_{n+1}))}$$
where $W_{p,q}(x_0,\ldots x_q) = \sum_{n=0}^{q-1} h(x_n,x_{n+1})$ on the space of sequences with $x_q=x_0+p$.  In the case of $(p,q)=(0,1)$ for the above fixed points, this gives $$\lambda+1/\lambda-2 = \frac{(1-g'(x))^2}{g'(x)},$$
so positive, and hence these fixed points of $\Phi$ are hyperbolic.

To complete the example, we show that the other branches of stable and unstable manifold of the two hyperbolic fixed points (shown in black and blue) can be chosen to be non-coincident.  


\begin{figure}[htbp] 
   \centering

   \begin{tikzpicture} [scale=0.9]
\draw [->](0,0)--(7,0)node[below]{$x$};
\draw [->](0.5,-1)--(0.5,1);
\draw [thick,red](0.5,0)..controls(2,-1)..(3.5,0);
\draw [thick,red](6.5,0)..controls(5,1)..(3.5,0);
\draw(1.25,-0.05)--(1.25,0.05);
\draw(2,-0.1)--(2,0.1)node[above]{$a_1$};
\draw(3.5,0.1)--(3.5,-0.1)node[below]{$b$};
\draw(3.5,-0.1)--(3.5,0.1)node[above]{$f(x)$};
\draw(5.75,-0.05)--(5.75,0.05);
\draw(5,0.1)--(5,-0.1)node[below]{$a_2$};
\draw(6.5,0.1)--(6.5,-0.1)node[below]{$1$};

\draw [->](8,0)--(15,0)node[below]{$t$};
\draw [->](8.5,-1)--(8.5,1);
\draw [thick,red](9.5,0)..controls(9.6,0)..(9.7,-0.1);
\draw [thick,red](9.7,-0.1)..controls(9.9,-0.3)..(10,0);
\draw [thick,red](10,0)..controls(10.1,0.3)..(10.3,0.1);
\draw [thick,red](10.3,0.1)..controls(10.4,0)..(10.5,0);

\draw [thick,red](13.5,0)..controls(13.4,0)..(13.3,-0.1);
\draw [thick,red](13.3,-0.1)..controls(13.1,-0.3)..(13,0);
\draw [thick,red](13,0)..controls(12.9,0.3)..(12.7,0.1);
\draw [thick,red](12.7,0.1)..controls(12.6,0)..(12.5,0);

\draw [thick,red](8.5,0)--(9.5,0);
\draw [thick,red](10.5,0)--(12.5,0);
\draw [thick,red](13.5,0)--(14.5,0);

\draw(9.25,-0.05)--(9.25,0.05);
\draw(10,-0.1)--(10,0.1)node[above]{$a_1$};
\draw(11.5,0.1)--(11.5,-0.1)node[below]{$b$};
\draw(13,0.1)--(13,-0.1)node[below]{$a_2$};
\draw(14.5,0.1)--(14.5,-0.1)node[below]{$1$};
\draw(13.75,-0.05)--(13.75,0.05);

\end{tikzpicture}

   \caption{(a) The graph of the function $f$; (b) the graph of the function $\beta$. }
   \label{fig:funfandbeta}
\end{figure}
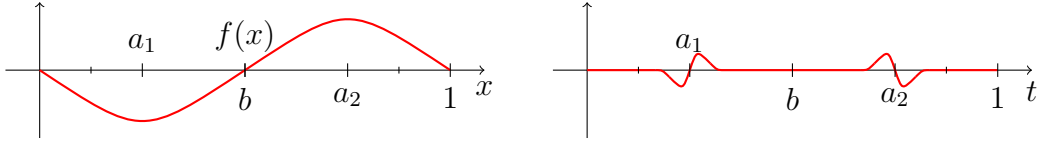

Let $g(x)=x+f(x)$. Assume $f'(x)\ge -1/2$ for all $x\in\R$ and $f$ has its minimum at $a_1\in (0,b)$ and maximum at $a_2\in (b,1)$ (see  Figure~\ref{fig:funfandbeta}(a)), that is, $f'(a_1)=f'(a_2)=0$, $f''(a_1)\ge 0$, and $f''(a_2)\le 0$. Furthermore, assume that
$$
f(x)=-\tfrac12 x\ \mbox{ for }\ x\in [0,a_1/2]\cup [(1+a_2)/2,1].
$$

Let $b_1=-f(a_1)>0$ and $b_2=-f(a_2)<0$. Then $\Phi$ has another two fixed points $(a_1,b_1)$ and $(a_2,b_2)$ other than $(0,0)$ and $(b,0)$. 

Construct a $C^1$ smooth function $\beta:\R\to\R$ which satisfies (see  Figure~\ref{fig:funfandbeta}(b))

(i)\ $\beta(t+1)=\beta(t)$, $\forall\,t\in\R$;

(ii)\ $\beta(a_1)=\beta(a_2)=0$, $\beta(t)\le 0$, $t\in [0,a_1]\cup [a_2,1]$, $\beta(t)\ge 0$, $t\in [a_1,a_2]$;

(iii)\ $|\beta(t)|<1/2 $ for $t\in U(a_1)\cup U(a_2)$ and $\beta(t)= 0$ for $t\not\in U(a_1)\cup U(a_2)$, where $U(a_1)$  and $U(a_2)$ are small neighborhoods of $a_1$ and $a_2$ respectively;

(iv)\ $\beta'(a_1)>2/b_1$ and $\beta'(a_2)<2/b_2$;

(v)\ $\int_b^x\beta(t)\mbox{d}t=0$ for $x\not\in U(a_1)\cup U(a_2)$;


Let $\hat{f}(x)=\int_b^x\beta(t)\mbox{d}t$, $\forall\,x\in\R$. Then $
\hat{f}(x+1)=\hat{f}(x)$, $\forall\,x\in\R$, in particular, $\hat{f}(x)=0$ for $x\not\in U(a_1)\cup U(a_2)$,
and $\hat{f}(a_1)<0$, $\hat{f}(a_2)>0$.

Now we perturb the original twist map $\Phi$ by letting 
$$
\tilde{f}(x)=f(x)+\hat{f}(x)\ \mbox{ and }\ g(x)=x+\tilde{f}(x).
$$
The new twist map, denoted by $\tilde{\Phi}$, has four fixed points: $(0,0), (b,0), (a_1,b_1')$ and $(a_2,b_2')$, where
$$
b_1'=-\tilde{f}(a_1)=b_1-\hat{f}(a_1)>b_1\ \mbox{ and }\ b_2'=-\tilde{f}(a_2)=b_2-\hat{f}(a_2)<b_2.
$$

In what follows we shall show that any point $(a_1,y)$ with $y\ge b_1$ is not on an invariant circle of $\tilde{\Phi}$. Assume the orbit of $\tilde{\Phi}$ starting from $(a_1,y)$ corresponds to a solution $\ubf=(u_n)$ of the difference equation
$$
-(x_n-g(x_{n-1}))+(x_{n+1}-g(x_n))g'(x_n)=0.
$$
Then $u_0=a_1$ and $y=u_1-g(u_0)=u_1-g(a_1)$.
Define a function
$$
W(z)=h(u_{-1},z)+h(z,u_1)=\frac{1}{2}(z-u_{-1}-\tilde{f}(u_{-1}))^2+\frac{1}{2}(u_1-z-\tilde{f}(z))^2.
$$
Then 
$$
W'(a_1)=0\ \mbox{ and }\ W''(a_1)=2-y\tilde{f}''(a_1)<0,
$$
due to the above condition (iv), implying that $\ubf=(u_n)$ is not a minimal configuration, i.e., the orbit starting from $(a_1,y)$ is not  minimal, and hence $(a_1, y)$ with $y\ge b_1$  does not lie on an invariant circle of the new twist map $\tilde{\Phi}$, since each orbit in an invariant circle of the twist map must be minimal~\cite{Ban,Mat0}. For the same reason, we deduce that $(a_2,y)$ with $y\le b_2$ is not on an invariant circle of $\tilde{\Phi}$ either.

As a consequence, neither the fixed point $(a_1,b_1')$ nor $(a_2, b_2')$ is  on an invariant circle of $\tilde{\Phi}$.

Note that $(a_1,3a_1/4)$ is on the unstable manifold of $(0,0)$ (the segment $y=3x/4, x\in[0,a_1]$ is in fact a part of the unstable manifold of $(0,0)$) and $3a_1/4>b_1$. Then the point $(a_1,3a_1/4)$ is not on an invariant circle, implying that the unstable manifold of $(0,0)$ and the stable manifold of $(b,0)$ do not coincide to form an invariant circle, and hence $F_d(0/1+)>0$ for the perturbed system, according to Theorem C. Similarly, we deduce that the stable  manifold of $(b,0)$ and the unstable manifold of $(1,0)$ do not coincide and hence $F_d(0/1-)>0$.

\section*{Appendix E:~Modification of $h$}
\label{sec:appE}
 \setcounter{theorem}{0}
 \setcounter{equation}{0}
 \setcounter{lemma}{0}
 \renewcommand{\theequation}{E.\arabic{equation}}
\renewcommand{\thelemma}{E.\arabic{lemma}}

Given $h \in C^2$ with $h(x+1,x'+1)=h(x,x')$, $h_{12}\le-c<0$, and integers $M<N$, we construct $\tilde{h} \in C^2$ with the same properties, equal to $h$ on $M\le x'-x\le N$ and with bounded second derivative.

We first modify for $x'-x>N$.  We can make a $C^1$ modification by defining $g(x)=h(x,x+N)$, $k(x)=h_2(x,x+N)-h_1(x,x+N)$ and using d'Alembert's solution of the wave equation with a constant source term, i.e.~for $x'-x\ge N$ let
$$\tilde{h}(x,x') = \frac12\left[g(x)+g(x'-N)+\int_x^{x'-N} k(\xi)\mbox{d}\xi + c(x'-x-N)^2\right].$$
In $x'-x> N$, $\tilde{h}$ is $C^2$ and has $\tilde{h}_{12} = -c$.
On $x'=x+N$, $\tilde{h}(x,x+N) = g(x)$, as desired.
It has
\begin{eqnarray}
\tilde{h}_1(x,x') &=& \tfrac12\left[g'(x) - k(x)\right]-c(x'-x-N), \nonumber\\
\tilde{h}_2(x,x') &=& \tfrac12\left[g'(x'-N)+k(x'-N)\right]+c(x'-x-N).\nonumber
\end{eqnarray}
On $x'=x+N$ these give
$$
\tilde{h}_1 = \tfrac12[h_1+h_2-h_2+h_1] = h_1\ \ \mbox{ and }\ \ \tilde{h}_2 = \tfrac12[h_1+h_2+h_2-h_1]=h_2,
 $$
 all evaluated at $(x,x+N)$, as desired.   Thus $\tilde{h}$ is a $C^1$ extension of $h$ from $x' -x\le N$ to $x'-x\ge N$.

To make a $C^2$ modification we have to be a bit more subtle.  Replace the term $c (x'-x-N)^2$ by
$$
4 \int_x^{x'-N} j(\xi)\min\{\xi-x,x'-N-\xi\}\mbox{d}\xi,
 $$
 where $j(x) =- h_{12}(x,x+N)$.  Then $\tilde{h}_{12}(x,x') = -j((x'-N+x)/2)$.  So now $\tilde{h}_{12}$ is a continuous extension of $h_{12}$.  Similarly we check that $\tilde{h}_{11}$ is a continuous extension of $h_{11}$:
$$\tilde{h}_1(x,x') = \frac12\left[g'(x)-k(x)-4\int_x^{(x'-N+x)/2} j(\xi) \mbox{d}\xi\right],$$
so
$$\tilde{h}_{11}(x,x') = \tfrac12[g''(x)-k'(x)-2j((x'-N+x)/2)+4j(x)] \to h_{11}(y,y+N)$$ as $(x,x') \to (y,y+N)$.
The same applies to $\tilde{h}_{22}$.

The analogous modification can be made for $x'<x+M$.

One could probably make smoother modifications, but the above suffices for present purposes.

\section*{Appendix F:~Representation for TOCs}
\label{app:TOCrep}
 \setcounter{theorem}{0}
 \setcounter{equation}{0}
 \setcounter{lemma}{0}
 \renewcommand{\theequation}{F.\arabic{equation}}
\renewcommand{\thelemma}{F.\arabic{lemma}}

Here we give a useful representation for TOCs and in particular, IOCs.

\begin{figure}[htbp] 
   \centering
   \begin{tikzpicture}[scale=0.7]
\draw (0,0)--(6,0)node [below] {$x_n$};
\draw (0,0)--(0,8.5);
\draw (0,0)--(6,6);
\draw (6,0)--(6,8.5);
\draw (6,6)--(0,6) node [anchor=east] {$x_{n+1}$};

\draw  (-0.3,1)--(1,1)--(1,3)--(3,3)--(3,7)--(6.2,7);
\draw [dashed](-0.3,2)--(2,2)--(2,5)--(5,5)--(5,8)--(6.2,8);

\begin{scope}[blue,thick]
\draw (0,2.5)..controls (0.5,2.7)..(1,3);
\draw (1,3)..controls (1.5,3.3)and(1.5,3.5)..(2,5);
\draw (2,5)..controls (2.5,6.5)and(2.8,6.9)..(3,7);
\draw (3,7)..controls (3.2,7.1)and(4.5,7.7)..(5,8);
\draw (5,8)..controls (5.5,8.3)..(6,8.5);
\begin{scope}[->]
\draw (0.4,2.67)--(0.5,2.7);
\draw (1.5,3.53)--(1.4,3.4);
\draw (2.26,5.75)--(2.33,5.9);
\draw (4.1,7.52)--(4,7.47);
\draw (5.5,8.28)--(5.6,8.35);
\end{scope}

\end{scope}

\end{tikzpicture}

   \caption{An IOC as the set of orbits of a circle homeomorphism, indicating the effect of the gradient flow.  This case has rotation number $1/2$ and two periodic orbits are shown.}
   \label{fig:circlehomeo}
\end{figure}
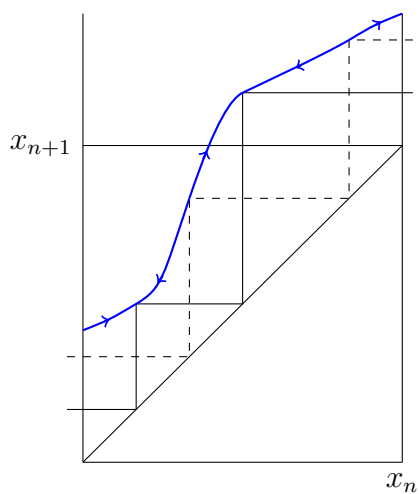

A TOC is fully specified by a $C^0$ strictly ordered curve $(x_0(s),x_1(s)) \in \R^2$, $s \in \R$, with the usual partial order on $\R^2$ and
\begin{equation}
(x_0(s+1),x_1(s+1)) = (x_0(s)+1,x_1(s)+1).
\label{eq:percond}
\end{equation}
This is because translation invariance lets us deduce $x_n(s)$ as follows: $x_{-1}(s) = x_0(s^-)$ where $s^-$ is determined by $x_1(s^-)=x_0(s)$ and $x_2(s) = x_1(s^+)$ where $s^+$ is determined by $x_0(s^+)=x_1(s)$.  The remaining $x_n(s)$ are determined by induction.  The parametrisation is not an intrinsic part of a TOC, however.  In particular, invariance of an IOC specified in the above way means the gradient flow induces a flow on the parameter $s$.
One can think of the strictly ordered curve as the graph of a lift of a degree-one circle-homeomorphism, and indeed the rule for determining a doubly infinite sequence $\xbf$ from a single point on the curve is precisely forward and backward iteration of the circle-homeomorphism.  The mean spacing for the sequences in a TOC is precisely the rotation number of the lifted circle-homeomorphism. See Figure~\ref{fig:circlehomeo}. 

\end{document}